\documentclass[twocolumn]{autart}    

\usepackage{graphicx}          
\usepackage{cite}
\usepackage{amsmath,amssymb,amsfonts,bm}
\usepackage{mathtools}
\usepackage{algorithmic}
\usepackage{graphicx}
\usepackage{textcomp}
\usepackage{dsfont}
\usepackage{epstopdf}        
\usepackage{enumerate}
\usepackage{lscape}
\usepackage{multicol}
\usepackage{multirow}
\usepackage{calligra}
\usepackage{booktabs}
\usepackage{mathtools}
\usepackage{framed}  
\usepackage{empheq}
\usepackage{graphics} 
\usepackage{epsfig} 
\usepackage{times} 
\usepackage{amsmath} 
\usepackage{amssymb}  
\usepackage{amsfonts}  
\usepackage{epstopdf}
\usepackage{enumerate}
\usepackage{graphicx} 
\usepackage{comment} 
\usepackage{apptools}
\usepackage{subcaption}
\usepackage{picins}
\usepackage{theoremref}
\usepackage{bbm}

\newtheorem{definition}{Definition}
\newtheorem{asmp}{Assumption}
\newtheorem{remark}{Remark}

\newcommand*{\QEDB}{\hfill\ensuremath{\square}}%
\def\BibTeX{{\rm B\kern-.05em{\sc i\kern-.025em b}\kern-.08em
    T\kern-.1667em\lower.7ex\hbox{E}\kern-.125emX}}

\begin{document}

\begin{frontmatter}

\vspace{-0.4cm}
\title{Hybrid Minimum-Seeking in Synergistic Lyapunov Functions: Robust Global Stabilization under Unknown Control Directions\thanksref{footnoteinfo}} 

\thanks[footnoteinfo]{This work was supported in part by the grants NSF ECCS CAREER 2305756 and AFOSR YIP: FA9550-22-1-0211. Corresponding Author: Jorge I. Poveda.}

\author[UCSD]{Mahmoud Abdelgalil}\ead{mabdelgalil@ucsd.edu},    
\author[UCSD]{Jorge I. Poveda}\ead{jipoveda@ucsd.edu},               
    
\address[UCSD]{Department of Electrical and Computer Engineering, University of California, San Diego, La Jolla, CA, USA.}


\begin{abstract}     
We study the problem of robust global stabilization in control-affine systems, focusing on dynamic uncertainties in the control directions \emph{and} the presence of topological obstructions that prevent the existence of smooth global control Lyapunov functions. Building on a recently developed Lie-bracket averaging result for hybrid dynamic inclusions presented in \cite{abdelgalil2023lie}, we propose a novel class of universal hybrid feedback laws that achieve robust global practical stability by identifying the minimum point of a set of appropriately chosen synergistic Lyapunov functions. As concrete applications of our results, we synthesize different hybrid high-frequency high-amplitude feedback laws for the solution of robust global stabilization problems on various types of manifolds under unknown control directions, as well as controllers for obstacle avoidance problems in vehicles characterized by kinematic models describing both holonomic and non-holonomic models. By leveraging Lie-bracket averaging for hybrid systems, we also show how the proposed hybrid minimum-seeking feedback laws can overcome lack of controllability during persistent (bounded) periods of time. Numerical simulation results are presented to illustrate the main results.
\end{abstract}

\end{frontmatter}

\section{Introduction}

\vspace{-0.2cm}
Ensuring robustness in the face of uncertainty is a core challenge in controller design for autonomous systems. This challenge becomes significantly more difficult if the uncertainty lies in the \textit{control direction}. For example, cyber-physical systems can experience significant damage when the control gain's sign unexpectedly changes or vanishes. Such a scenario may occur without malicious interference, such as from internal software failures, or due to external spoofing attacks by adversarial agents \cite{pasqualetti2013attack}. Furthermore, the control direction may be time-varying and intermittently zero (indicating a lack of controllability) for bounded periods. In such cases, traditional identification or parameter estimation techniques are often inadequate due to insufficient excitation in the system's trajectories. As a result, designing stabilizing feedback laws capable of addressing uncertainty in control directions and persistent periods of uncontrollability is essential for ensuring the resilience of autonomous systems operating in complex and dynamic environments.

Referred to as the problem of \emph{``stabilization with an unknown sign of the high-frequency gain''}, or \emph{``stabilization under unknown control directions''}, designing stabilizing feedback laws in the presence of uncertainty on control directions has a long history in the literature of robust adaptive control that dates back to the 1970s \cite{morse1979open}. The first solution was proposed by Nussbaum \cite{nussbaum1983some}, which became the common approach to solving this problem. However, it is well-known, e.g. see \cite{AlexTAC}, that the Nussbaum approach can lead to poor performance if the control sign is dynamic. An alternative \emph{model-free} approach for solving this difficult problem was introduced by Scheinker and Krstic in \cite{AlexTAC,scheinker2017model}. The approach introduced in \cite{AlexTAC} is based on \emph{seeking} the minimum of a suitably constructed \emph{control Lyapunov-like function} $V$, using a high-frequency time-varying oscillatory feedback law that relies only on real-time evaluations of this function. As shown in \cite{AlexTAC,scheinker2017model}, this model-free approach is impervious to persistent changes in the sign of the control gain, which makes it an ideal candidate for the design of stabilizing feedback laws operating in uncertain and adversarial environments. Since its conception in \cite{AlexTAC}, the minimum-seeking approach has found numerous applications, including in the context of 2-D Vehicle Control \cite{scheinker2017model}, output regulation problems in nonlinear systems \cite{wang2024extremum}, etc, see \cite[Sec. 6]{scheinker2024100} for a recent survey on this subject.

\vspace{-0.3cm}
A core theoretical development that enabled the results of \cite{AlexTAC,scheinker2017model} was the emergence of \emph{Lie-bracket averaging} for ODEs as a powerful framework for the analysis and design of model-free control and optimization algorithms \cite{Durr:13,grushkovskaya2018class,suttner2017exponential,suttner2022extremum}. Nevertheless, since the feedback laws analyzed in \cite{AlexTAC, scheinker2017model} are continuous and derive their stability properties from the stability of their Lie-bracket averaged systems (which are also continuous), their performance is inherently constrained by the standard limitations associated with continuous feedback laws \cite{sontag_discontinuous}. In particular, when the system operates in topological spaces that are not contractible to a point, as in, e.g., smooth compact manifolds, robust global stabilization of a desired point is not possible via continuous feedback \cite{jongeneel2023topological}. Such types of topological obstructions to robust global stabilization have been shown to emerge in different applications, including obstacle avoidance problems \cite{Poveda:20TAC,TrackingManifold}, the attitude stabilization of rigid bodies\cite{SynchronizationMeyhew,mayhew2011synergistic}, stabilization problems in $\mathbb{S}^1$, synchronization problems \cite{SynchronizationMeyhew}, etc. When the sign in the control direction is static and \emph{known} a priori, such stabilization problems can be tackled via synergistic hybrid control \cite{Mayhew10Thesis,TrackingManifold,RSanfeliceBook}. However, most existing results on synergistic hybrid control cannot cope with uncertainty in the control directions, particularly when such uncertainty is dynamic. Indeed, even the state-of-the-art results on synergistic hybrid feedback for uncertain systems \cite{casau2024robust} are inapplicable when the control gain sign is uncertain.

\vspace{-0.1cm}
Motivated by the above challenges, the main contribution of this paper is to provide a solution for \textit{robust global stabilization problems under unknown control directions in spaces that are not contractible to a point}. Specifically, the following are the  main contributions of this paper:

\vspace{-0.2cm}\noindent
1) We propose a novel class of hybrid and oscillatory feedback laws that leverage the existence of a collection of suitable \emph{local} strong control Lyapunov functions (SCLFs) available to the plant. These local SCLFs are then used together to control the system when a global SCLF is unavailable. Building on existing results in model-based hybrid control \cite{Mayhew10Thesis,mayhew2011synergistic,TrackingManifold} and using a recently introduced Lie-bracket averaging result for hybrid inclusions \cite{abdelgalil2023lie}, we establish \emph{semi-global practical asymptotic stability} for systems with unbounded operational sets, and \emph{global practical asymptotic stability} for systems on smooth compact manifolds, even with unknown control directions that vanish intermittently. 

\vspace{-0.2cm}\noindent
2) Subsequently, we apply our results to various control problems with unknown control directions where a robustly globally stabilizing continuous control law does not exist even if the control direction is known, including stabilization problems on $\mathbb{S}^1$, $\mathbb{S}^2$, $\text{SO}(3)$, and stabilization problems on the plane with obstacle avoidance for both holonomic and non-holonomic vehicles. 

\vspace{-0.2cm}\noindent
3) Finally, we present different numerical examples to demonstrate the performance of the proposed controllers under different types of uncertainty in the control directions. Additionally, we compare their effectiveness against both non-hybrid approaches and hybrid non-adaptive methods. In all applications, we show how our results open the door for the design of ``model-free'' controllers based on existing well-posed model-based hybrid algorithms \cite{Mayhew,RSanfeliceBook,TrackingManifold,mayhew2010hybrid}, thus demonstrating how the synergistic use of hybrid control \cite{RSanfeliceBook} and Lie-bracket averaging \cite{abdelgalil2023lie} can simultaneously overcome uncertainties in control directions \emph{and} the topological obstructions to global stabilization.

To the best of our knowledge, the results presented in this paper provide the first theoretical link between the minimum-seeking approach introduced in \cite{scheinker2017model}, and the setting of synergistic hybrid control \cite{Mayhew10Thesis,RSanfeliceBook}. 

\vspace{-0.1cm}
The remainder of this manuscript is organized as follows. In Section \ref{sec:prelims}, we present the preliminaries. Section \ref{sec:SCLF} presents the main problem formulation and main theoretical result. Section \ref{sec:Synergistic_Lyapunov_Functions} focuses on robust global stabilization problems on smooth manifolds. The proofs are presented in Section \ref{sec:Proofs}, followed by the conclusions and future work in Section \ref{sec:conclusions}.

\vspace{-0.3cm}
\section{Preliminaries}\label{sec:prelims}

\vspace{-0.2cm}
%
\subsection{Notation}

\vspace{-0.2cm}
We use $\langle x, y\rangle = x^\top y$, to denote the inner product between any two vectors $x,y\in\mathbb{R}^n$. Given a compact set $\mathcal{A}\subset\mathbb{R}^n$ and  $x\in\mathbb{R}^n$, we use $|x|_{\mathcal{A}}:=\min_{\tilde{x}\in\mathcal{A}}\|x-\tilde{x}\|_2$. A set-valued mapping $M:\mathbb{R}^p\rightrightarrows\mathbb{R}^n$ is outer semicontinuous (OSC) at $z$ if for each sequence $\{z_i,s_i\}\to(z,s)\in\mathbb{R}^p\times\mathbb{R}^n$ satisfying $s_i\in M(z_i)$ for all $i\in\mathbb{Z}_{\geq0}$, we have $s\in M(z)$. A mapping $M$ is locally bounded (LB) at $z$ if there exists an open neighborhood $N_z\subset\mathbb{R}^p$ of $z$ such that  $M(N_z)$ is bounded. The mapping $M$ is OSC and LB relative to a set $K\subset\mathbb{R}^p$ if $M$ is OSC for all $z\in K$ and $M(K):=\cup_{z\in K}M(x)$ is bounded. A function $\beta:\mathbb{R}_{\geq0}\times\mathbb{R}_{\geq0}\to\mathbb{R}_{\geq0}$ is of class $\mathcal{K}\mathcal{L}$ if it is nondecreasing in its first argument, nonincreasing in its second argument, $\lim_{r\to0^+}\beta(r,s)=0$ for each $s\in\mathbb{R}_{\geq0}$, and  $\lim_{s\to\infty}\beta(r,s)=0$ for each $r\in\mathbb{R}_{\geq0}$. For two (or more) vectors $u,v \in \mathbb{R}^{n}$, we write $(u,v)=[u^{\top},v^{\top}]^{\top}$. If $x,y\in\mathbb{R}^3$, we use $[x]_{\times}$ to denote the the skew-symmetric matrix associated to $x$ and defined such that $[x]_\times y = x\times y$, where $\times$ indicates the cross product of vectors in $\mathbb{R}^3$. The notation $\otimes$ denotes the Kronecker product. Let $x=(x_1,x_2,\ldots,x_n)\in\mathbb{R}^n$. Given a Lipschitz continuous function $f:\mathbb{R}^n\to\mathbb{R}^m$, we use $\partial_{x_i}f$ to denote the generalized Jacobian \cite{clarke2008nonsmooth} of $f$ with respect to the variable $x_i$. A map $f$ is said to be of class $\mathcal{C}^k$ if it is $k$-times continuously differentiable with the $k$th-derivative being locally Lipschitz continuous. If $C\subset\mathbb{R}^n$, the notation $T_x C$ denotes the \emph{tangent cone} \cite{clarke2008nonsmooth} of $C$ at the point $x$. When $C$ is an embedded Euclidean submanifold, $T_x C$ denotes the tangent space at $x$, which is isomorphic to an affine linear subspace of $\mathbb{R}^n$. 

\vspace{-0.2cm}
\subsection{Hybrid Dynamical Systems}

\vspace{-0.2cm}
In this paper, our models are given by hybrid dynamical systems (HDS), as studied in \cite{bookHDS}. Such systems are characterized by the following inclusions:
\begin{subequations}\label{eq:HDS0}
\begin{align}
\mathcal{H}:~~~\begin{cases}
    ~~x\in C, & \dot{x}\hphantom{^+}\in F(x)\\
    ~~x\in D, & x^+\in G(x),
\end{cases}
\end{align}
\end{subequations}
where $F:\mathbb{R}^n\rightrightarrows\mathbb{R}^n$ is called the flow map, $G:\mathbb{R}^n\rightrightarrows\mathbb{R}^n$ is called the jump map, $C\subset\mathbb{R}^n$ is called the flow set, and $D\subset\mathbb{R}^n$ is called the jump set. We use $\mathcal{H}=(C,F,D,G)$ to denote the \emph{data} of the HDS $\mathcal{H}$. Purely continuous-time systems can be modeled as \eqref{eq:HDS0} by taking $D=\{\emptyset\}$. Similarly, purely discrete-time systems can be modeled as \eqref{eq:HDS0} by taking $C=\{\emptyset\}$.  In this paper, we work with well-posed HDS that satisfy the following assumption \cite[Assumption 6.5]{bookHDS}.
\begin{asmp}\label{asmp:basic}[Hybrid Basic Conditions]
The sets $C,D$ are closed. The set-valued mapping $F$ is OSC, LB, and for each $x\in C$ the set $F(x)$ is convex and nonempty. The set-valued mapping $G$ is OSC, LB, and for each $x\in D$ the set $G(x)$ is nonempty. \QEDB 
\end{asmp}
Henceforth, all hybrid systems in this manuscript are assumed to satisfy the Hybrid Basic Conditions.

Solutions to system \eqref{eq:HDS0} are parameterized by a continuous-time index $t\in\mathbb{R}_{\geq0}$, which increases continuously during flows, and a discrete-time index $j\in\mathbb{Z}_{\geq0}$, which increases by one during jumps. Therefore, solutions to \eqref{eq:HDS0} are defined on \emph{hybrid time domains} (HTDs). A set $E\subset\mathbb{R}_{\geq0}\times\mathbb{Z}_{\geq0}$ is called a \textsl{compact} HTD if $E=\cup_{j=0}^{J-1}([t_j,t_{j+1}],j)$ for some finite sequence of times $0=t_0\leq t_1\ldots\leq t_{J}$. The set $E$ is a HTD if for all $(T,J)\in E$, $E\cap([0,T]\times\{0,\ldots,J\})$ is a compact HTD. The following definition formalizes the notion of solution to HDS of the form \eqref{eq:HDS0}.
\begin{definition}\normalfont \label{definitionsolutions1}
A hybrid arc $x$ is a function defined on a HTD. In particular, $x:\text{dom}(x)\to \mathbb{R}^n$ is such that $x(\cdot, j)$ is locally absolutely continuous for each $j$ such that the interval $I_j:=\{t:(t,j)\in \text{dom}(x)\}$ has a nonempty interior. A hybrid arc $x:\text{dom}(x)\to \mathbb{R}^n$ is a solution $x$ to the HDS \eqref{eq:HDS0} if $x(0, 0)\in C\cup D$, and: 1) For all $j\in\mathbb{Z}_{\geq0}$ such that $I_j$ has nonempty interior: $x(t,j)\in C$ for all $t\in I_j$, and $\dot{x}(t,j)\in F(x(t,j))$ for almost all $t\in I_j$; 2) For all $(t,j)\in\text{dom}(x)$ such that $(t,j+1)\in \text{dom}(x)$: $x(t,j)\in D$ and $x(t,j+1)\in G(x(t,j))$. A solution $x$ is said to be \emph{maximal} if it cannot be further extended, and it is said to be \emph{complete} if $\text{length}~\text{dom}(x)=\infty$. \QEDB
\end{definition}


%

In this paper, we use the following standard stability notion:
\begin{definition}\normalfont \label{definitionstablity1}
    A compact set $\mathcal{A}\subset\mathbb{R}^n$ is said to be \emph{uniformly globally asymptotically stable (UGAS)} for the HDS \eqref{eq:HDS0} if there exists $\beta\in\mathcal{KL}$ such that each solution $x$ to \eqref{eq:HDS0} satisfies $|x(t,j)|_{\mathcal{A}}\leq \beta(|x(0,0)|_{\mathcal{A}},t+j)$, for all $(t,j)\in\text{dom}(x)$. \QEDB 
\end{definition}
%
%
%
In some cases, we consider HDS that depend on a small tunable parameter $\varepsilon\in\mathbb{R}_{>0}$, given by
\begin{subequations}\label{eq:ms_orig_hybrid_sys}
    \begin{align}
    \mathcal{H}_\varepsilon:
        \begin{cases}    
        C & \dot{x}\hphantom{^+}\in F_\varepsilon(x) \\
        D & x^+\in G_\varepsilon(x).
        \end{cases}
    \end{align}
\end{subequations}
For system \eqref{eq:ms_orig_hybrid_sys}, we use the following stability notion, which is standard in the literature \cite{Wang:12_Automatica,scheinker2017model}.
\begin{definition}\normalfont \thlabel{defn:SPUAS}
    For the HDS (\ref{eq:ms_orig_hybrid_sys}), a compact set $\mathcal{A}\subset\mathbb{R}^{n}$ is said to be \emph{Semi-Globally Practically Asymptotically Stable} (SGpAS) as $\varepsilon\rightarrow 0^+$ if there exists $\beta\in\mathcal{KL}$ such that for each compact set $K\subset C\cup D$ and for each $\nu>0$, there exists $\varepsilon^*>0$ such that for all $\varepsilon\in(0,\varepsilon^*]$, all solutions of \eqref{eq:ms_orig_hybrid_sys} with $x(0,0)\in K$ satisfy $|x(t,j)|_{\mathcal{A}}\leq \beta(|x(0,0)|_{\mathcal{A}},t+j)+ \nu$, for all $(t,j)\in\text{dom}(x)$. \QEDB 
\end{definition}
If there are several parameters that need to be tuned sequentially, we use the notation $(\varepsilon_1,\varepsilon_2,\dots,\varepsilon_\ell)\rightarrow 0^+$ to encode the sequence of parameter tuning, i.e. $\varepsilon_1$ is tuned, then $\varepsilon_2$ after fixing $\varepsilon_1$, and so forth. Finally, note that when $C\cup D$ is a compact set, SGpAS reduces to \emph{Uniform Global Practical Asymptotic Stability} (UGpAS). 

\section{Strong $\nabla V$-Stabilizability of Control-Affine HDS}\label{sec:SCLF}
In this section, we present the model of the systems under study, and a general result on model-free stabilization of control-affine HDS \emph{with unknown control directions}. 
\subsection{Model}
Consider the open-loop HDS:
\begin{subequations}\label{eq:org_HDS}
    \begin{align}
        \begin{cases}~
            x\in\mathcal{X}_C,&
            \dot{x}\hphantom{^+} = f(x,\theta,u)
        \\
        ~x\in \mathcal{X}_D, & x^+\in G_{\mathcal{X}}(x)
        \end{cases}
    \end{align}
where $x\in\mathbb{R}^n$ is the main state, $u := (u_1,u_2,\dots,u_{r})\in \mathbb{R}^{r}$ are the control inputs, $\mathcal{X}_C \subset \mathbb{R}^n$ is the flow set, $\mathcal{X}_D \subset\mathbb{R}^n$ is the jump set, $G_{\mathcal{X}}:\mathbb{R}^n\rightrightarrows\mathbb{R}^n$ is the jump map, and $f:\mathbb{R}^n\times\mathbb{R}^m\times\mathbb{R}^r$ is the flow map, given by
\begin{align}\label{flowmap00}
        f(x,\theta,u) &:= f_0(x,\theta) + \sum_{i=1}^r f_i(x,\theta) u_i.
    \end{align}
\end{subequations}
In \eqref{flowmap00}, $\theta\in\Theta\subset\mathbb{R}^m$ represents a vector of potentially time-varying \emph{unknown parameters} that models the unknown control directions. The class of systems \eqref{eq:org_HDS} generalizes the class of \emph{control-affine} systems commonly studied in the literature on continuous-time systems modeled as ODEs \cite{scheinker2017model} by allowing for jumps in the solutions. We can also consider the class of systems \eqref{flowmap00} as an open-loop HDS in the sense of \cite{sanfelice2013existence}, where the flow set, jump set, and jump map have been designed, and it remains to design the continuous feedback laws $u_i$. We assume that \eqref{eq:org_HDS} satisfies the Hybrid Basic Conditions.

The following regularity condition on \eqref{eq:org_HDS} is needed to guarantee existence of maximal solutions and to rule out pathological behaviors, such as  Zeno solutions. 
\begin{asmp}\label{asmp:regularity_org_HDS} For all $i\in\{1,2,\dots,r\}$ and all $(x,\theta)\in\mathcal{X}_C\times\Theta$, the following holds
\begin{enumerate}
    \item $f_0$ and $f_i$ are $\mathcal{C}^0$, and $f_i(\cdot,\theta)$ is $\mathcal{C}^1$.
    \item $f_0(x,\theta)\in T_{x}\mathcal{X}_C$, and $f_i(x,\theta)\in T_{x}\mathcal{X}_C$.
    \item $G(\mathcal{X}_D)\subset \mathcal{X}_C$, and $G(\mathcal{X}_D)\cap \mathcal{X}_D = \emptyset$.
\end{enumerate}
\end{asmp}

The vector of unknown parameters $\theta$ may be constant or time-varying. In the latter case, we model $\theta$ as the solution of an \emph{exogenous} HDS of the form
\begin{align}\label{eq:exogenous_signal_dynamics}
    \begin{cases}
        ~~\theta\in \Theta_C&
        \dot{\theta}\hphantom{^+} \in F_e(\theta)
    \\
    ~~\theta\in \Theta_D & \theta^+ \in G_e(\theta),
    \end{cases}
\end{align}
where $\Theta_C,\Theta_D\subset \Theta$. We remark that the case of static uncertainty $\dot{\theta}=0$ is trivially included in \eqref{eq:exogenous_signal_dynamics} by taking $\Theta_C=\{\Theta\}$, $\Theta_D=\emptyset$, and $F_e(\theta)=\{0\}$. As usual, we require that the HDS \eqref{eq:exogenous_signal_dynamics} satisfies the Hybrid Basic Conditions of Assumption \ref{asmp:basic}. In addition, to guarantee existence of solutions to \eqref{eq:exogenous_signal_dynamics} from any initial condition in $\Theta_C\cup \Theta_D$, we impose the following regularity assumption.
\begin{asmp}\label{asmp:exogeneous_signal_regularity} a) For each $\vartheta\in \Theta_C$ there exists a neighborhood $U$ of $\vartheta$ such that, for all $\theta\in \Theta_C\cap U$ we have $F_e(\theta)\cap T_{\theta}\Theta_C\neq\emptyset$; b)
$G_e(\Theta_D)\subset \Theta_C$; d) and $G_e(\Theta_D)\cap \Theta_D = \emptyset$.
\end{asmp}
Under Assumption \ref{asmp:exogeneous_signal_regularity}, the HDS defined by \eqref{eq:exogenous_signal_dynamics} is sufficiently general to model a variety of complex behaviors. For example, the HDS \eqref{asmp:exogeneous_signal_regularity} can model various classes of switching signals \cite[Section 2.4]{bookHDS}, including those that have dwell-time bounds, which can be generated using hybrid automata \cite{PoTe17Auto,liu2025impulsive}. While item (a) in  Assumption \ref{asmp:exogeneous_signal_regularity} rules out consecutive jumps of $\theta$, this assumption can be relaxed to consider systems that satisfy average dwell-time bounds. However, since they rarely occur in practice, we rule out consecutive jumps of the control direction.

\begin{remark}\normalfont
While $\theta$ is treated in this paper as an uncertain parameter, the regularity properties on the dynamics of $\theta$ are necessary so as to exclude pathological behaviors, e.g. Zeno solutions, purely discrete solutions, etc. In particular, by \cite[Proposition 2.34]{sanfelice2021hybrid}, Assumption \ref{asmp:exogeneous_signal_regularity} guarantees that, due to the compactness of $\Theta$, there exists $t_\circ>0$ such that any solution $\theta$ to \eqref{eq:exogenous_signal_dynamics} satisfies $\overline{t}_{j}-\underline{t}_{j}\geq t_\circ > 0$ for all $j$ in the domain of $\theta$, where $\overline{t}_{j}:=\sup\{t\in\mathbb{R}_{\geq 0}~:~(t,j)\in\text{dom}(\theta)\}$, $\underline{t}_{j}:=\,\inf\,\{t\in\mathbb{R}_{\geq 0}~:~(t,j)\in\text{dom}(\theta)\}$. \QEDB 
\end{remark}
By interconnecting systems \eqref{eq:org_HDS} and \eqref{eq:exogenous_signal_dynamics}, we obtain the following open-loop HDS with unknown control directions:
\begin{subequations}\label{eq:combined_org_HDS}
\begin{align}
    \mathcal{H}:~\begin{cases}
        ~~\xi\in C, & \dot{\xi}\hphantom{^+}\in F(\xi,u)\\
        ~~\xi\in D, & \xi^+\in G(\xi),
    \end{cases}
\end{align}
where $\xi=(x,\theta)$, and where the sets $C$ and $D$ are given by
\begin{align}\label{constructionflowjumpsets}
        C &= \mathcal{X}_C\times \Theta_C,\\
        D &= (\mathcal{X}_C\times\Theta_D)\cup  (\mathcal{X}_D\times\Theta_C)\cup  (\mathcal{X}_D\times\Theta_D),
\end{align}
the flow map $F$ is given by
\begin{align}
        F(\xi,u)&= \{f(x,\theta,u)\}\times F_e(\theta),
\end{align}
and the jump map $G$ is given by
\begin{align}\label{jumpmapGd}
        G(\xi)=\begin{cases}
            G_{\mathcal{X}}(x)\times G_e(\theta) & (x,\theta)\in \mathcal{X}_D\times \Theta_D\\
            \{x\}\times G_e(\theta) & (x,\theta)\in \mathcal{X}_C\times \Theta_D\\
            G_{\mathcal{X}}(x)\times\{\theta\} & (x,\theta)\in \mathcal{X}_D\times \Theta_C.
        \end{cases}
    \end{align}
\end{subequations}
In words, the jump map \eqref{jumpmapGd} allows jumps in the system whenever $x$, or $\theta$, are in their respective jump sets $\mathcal{X}_D$ or $\Theta_D$, respectively. Note that, by construction, $\mathcal{H}$ satisfies Assumption \ref{asmp:basic} whenever the individual elements of \eqref{eq:org_HDS} and \eqref{eq:exogenous_signal_dynamics} satisfy the Hybrid Basic Conditions. In addition, the following Lemma is straightforward to verify. 
\begin{lem}\label{lem:solution_properties_lemma_aux}
    Suppose that Assumptions \ref{asmp:regularity_org_HDS} and \ref{asmp:exogeneous_signal_regularity} hold. Then, $G(D)\subset C$, and $G(D)\cap D = \emptyset$.
\end{lem}
\subsection{SCLFs for HDS with Unknown Control Directions}
To study the stabilization problem of system \eqref{eq:combined_org_HDS}, and to simplify some expressions, in the sequel we let $\mathcal{X}=\mathcal{X}_C\cup\mathcal{X}_D$, $\Theta=\Theta_C\cup\Theta_D$, 
%
%
and $\pi_x:\mathbb{R}^n\times\mathbb{R}^m \ni \xi=(x,\theta)\mapsto x\in \mathbb{R}^n$ be the canonical projection onto the first state. 
The following definition is inspired by the notion of \emph{``strong $L_g V$ stabilizability"}, introduced in \cite{scheinker2017model,AlexTAC} for input-affine ODEs.
\begin{definition}\normalfont \label{defn:strong_stabilizability}
    Let $\mathcal{A}\subset \mathcal{X}$ be a compact set. The $\mathcal{C}^1$ function $V:\mathbb{R}^n\rightarrow\mathbb{R}_{\geq 0}$ is said to be a \emph{Strong Control Lyapunov Function (SCLF) candidate with respect to $\mathcal{A}$ for $\mathcal{H}$} if there exists $\gamma>0$ and class $\mathcal{K}_\infty$-functions $\alpha_1$, $\alpha_2$ such that:
    \begin{enumerate}[(a)]
    \item For all $x\in \mathcal{X}$, we have:
    \begin{subequations}\label{Lyapunovcompactnotation}
    \begin{align}
     \alpha_1(|x|_{\mathcal{A}})\leq \,V(x) \leq \alpha_2(|x|_{\mathcal{A}}).
    \end{align}
    \item For all $\xi \in C$, we have:
    \begin{align}\label{flowconditionCLF}
    \dot{V}(\xi):= \langle \nabla V(x),\bar{f}(\xi)\rangle\leq 0,
    \end{align}
    where $\bar{f}$ is given by
\begin{align}\label{eq:avg_control_affine_sys}
    \bar{f}(\xi)&:= f_0(x,\theta)- \gamma \sum_{i=1}^r\langle \nabla V(x), f_i(x,\theta)\rangle f_i(x,\theta).
\end{align} 

\vspace{-0.3cm}
    \item For all $\xi\in D$, we have:
    \begin{align}
\Delta V(\xi):=\max_{g\in G(\xi)}V(\pi_x(g))-V(\pi_x(\xi))\leq 0.
      \end{align}
    \end{subequations}
    \end{enumerate}
\end{definition}
\begin{remark}\normalfont
Definition \ref{defn:strong_stabilizability} aims to generalize to HDS the \emph{``strong $L_g V$ stabilizability"} property studied in \cite{scheinker2017model,AlexTAC} for ODEs. Indeed, condition \eqref{flowconditionCLF} can be written in a more explicit form as
\begin{align}\label{flowweak}
    \langle\nabla V(x), f_0(x,\theta)\rangle - \gamma \sum_{i=1}^r\langle \nabla V(x),f_i(x,\theta)\rangle^2 \leq 0,
\end{align}
for all $(x,\theta)\in C$. If, in addition, for all $(x,\theta)\in C$ such that $|x|_{\mathcal{A}}=\epsilon > 0$, the function $V$ satisfies the \emph{``strong small control property"} \cite[Sec. 3.1]{scheinker2017model}:
\begin{align*}
    \lim_{\epsilon\rightarrow 0} \max_{\langle\nabla V(x), f_0(x,\theta)\rangle > 0}\frac{\langle\nabla V(x), f_0(x,\theta)\rangle}{\sum_{i=1}^r\langle \nabla V(x),f_i(x,\theta)\rangle^2} < +\infty,
\end{align*}
then, as shown in \cite{AlexTAC}, it is possible to construct another SCLF candidate for which the non-increase inequality in \eqref{flowweak} is strengthened to a strict decrease:
\begin{align*}
    \langle\nabla V(x), f_0(x,\theta)\rangle - \gamma \sum_{i=1}^r\langle \nabla V(x),f_i(x,\theta)\rangle^2 \leq -\alpha(|x|_{\mathcal{A}}),
\end{align*}
for some positive definite function $\alpha$. However, the new function will involve the expressions $\langle \nabla V(x),f_i(x,\theta)\rangle$, which necessitates the measurement of the unknown parameter $\theta$. Thus, we do not insist upon the strong small control property for $V$ so as to allow for some generality on how the unknown parameter $\theta$ affects the dynamics. \QEDB 
\end{remark}

Given a function $V$ that is a SCLF \emph{candidate} with respect to $\mathcal{A}$ for $\mathcal{H}$, we introduce the following auxiliary HDS:
\begin{align}\label{eq:org_HDS_V}
    \mathcal{H}_V:\begin{cases}
        ~~\xi\in C,&
        \dot{\xi}\hphantom{^+} \in \bar{F}(\xi)
    \\
    ~~\xi\in D, & \xi^+\in G(\xi)
    \end{cases}
\end{align}
wherein the flow map $\bar{F}$ is defined by
\begin{align*}
    \bar{F}(\xi)&:= \{\bar{f}(\xi)\}\times F_e(\theta),
\end{align*}
the map $\bar{f}$ is given by \eqref{eq:avg_control_affine_sys}, and the remaining data of $\mathcal{H}_V$ coincide with the data of $\mathcal{H}$. Note that the HDS $\mathcal{H}_V$ in \eqref{eq:org_HDS_V} is obtained by closing the loop for the HDS $\mathcal{H}$ in \eqref{eq:combined_org_HDS} with the following \emph{ideal} feedback law:
\begin{align}\label{eq:avg_feedback_law}
        u_i(\xi) = -\gamma \langle\nabla V(x),f_i(\xi)\rangle.
    \end{align}
By construction and inequalities \eqref{Lyapunovcompactnotation}, the compact set $\mathcal{A}\times\Theta$ is Lyapunov stable for $\mathcal{H}_V$ \cite[Theorem 3.18]{sanfelice2021hybrid}. However, implementing the ideal feedback law \eqref{eq:avg_feedback_law}, which is ubiquitous in the literature, requires real-time measurements of the state $x$ and the control direction $\theta$, as well as complete knowledge of $V$ and $f_i$. This level of information, however, is not accessible to the control law in \eqref{eq:combined_org_HDS}, leading to the following observation:

\vspace{0.0cm}\noindent 
\textbf{Fact:} For the HDS \eqref{eq:combined_org_HDS} with unknown control directions, the ideal feedback law \eqref{eq:avg_feedback_law} is not implementable. \QEDB 

Even though the SCLF-based control law $u_i$ is not suitable for implementations, our main goal is to use information about the existence of an SCLF to design an implementable controller able to \emph{emulate} $u_i$ and stabilize system \eqref{eq:combined_org_HDS}. However, since SCLF candidates only guarantee stability for the ideal closed-loop system $\mathcal{H}_V$, additional structure might be needed to obtain asymptotic stability. To capture such structure, we introduce the notion of \emph{strong $\nabla V$-stabilizability}, which is also instrumental to conclude when $V$ is an actual SCLF (not just a ``candidate'') with respect to $\mathcal{A}$ for $\mathcal{H}$.

\begin{definition}\normalfont \label{defn:strong_stabilizability_2}
    The HDS $\mathcal{H}$ is said to be \emph{strongly $\nabla V$-stabilizable} if there exists a SCLF candidate $V$ with respect to $\mathcal{A}$ for $\mathcal{H}$ such that $\mathcal{A}\times\Theta$ is UGAS for the ideal closed-loop system $\mathcal{H}_V$. When $\mathcal{H}$ is strongly $\nabla V$-stabilizable, $V$ is said to be a \emph{Strong Control Lyapunov Function (SCLF) with respect to $\mathcal{A}$\, for $\mathcal{H}$}. 
\end{definition}
\begin{remark}\normalfont
The topology of the problem under consideration may preclude the existence of a SCLF $V$. Indeed, this is the case for stabilization problems defined on smooth compact manifolds, where (robust) global stabilization via continuous feedback is not possible \cite{jongeneel2023topological}. However, we will show in Section \ref{sec:Synergistic_Lyapunov_Functions} that this issue can be overcome by considering a \emph{collection of functions}  $\{V_q\}_{q\in\mathcal{Q}}$, indexed by a discrete state $q$ that acts as a logic mode to be selected in real time by the controller. In this case, the state $x$ will be decomposed into $x=(p,q)$, where $p$ belongs to the smooth manifold and $q$ belongs to a finite discrete set $\mathcal{Q}$. It will be shown that the conditions of Definition \ref{defn:strong_stabilizability_2} are satisfied by the function $V(x)=V_q(p)$, which can be used for the purpose of real-time control using hybrid feedback. \QEDB 
\end{remark}

\begin{remark}\normalfont
    Establishing that the open-loop HDS $\mathcal{H}$ is strongly $\nabla V$-stabilizable is an application-dependent task. A sufficient condition is to establish the existence of a positive definite function $\alpha$ such that the following inequalities hold along the solutions of the ideal closed-loop HDS $\mathcal{H}_V$:
\begin{subequations}\label{eq:SCLF_suff_condition_PD}
        \begin{align}
            \dot{V}(\xi)&\leq -\alpha(|\pi_x(\xi)|_{\mathcal{A}}), & \forall \xi&\in C,\\
            \Delta V(\xi)&\leq -\alpha(|\pi_x(\xi)|_{\mathcal{A}}), & \forall \xi&\in D.
        \end{align}
    \end{subequations}
    However, condition \eqref{eq:SCLF_suff_condition_PD} place significant restrictions on the unknown parameter $\theta$. Thus, in lieu of \eqref{eq:SCLF_suff_condition_PD}, we only require that $V$ is a \emph{weak Lyapunov function} \cite{sanfelice2021hybrid} for $\mathcal{H}_V$ as described by the inequalities \eqref{Lyapunovcompactnotation}.Therefore, the weak decrease of $V$ may be used in conjunction with the properties of the solutions of $\mathcal{H}_V$ to certify that $\mathcal{A}\times(\Theta_C\cup\Theta_D)$ is UGAS for $\mathcal{H}_V$. Among others, this allows to model the practical situation wherein \emph{temporary loss of control} inhibits the ability of feedback to induce strict decrease of $V$ during flows and/or jumps. We further expand on this case in Section \ref{sec:Synergistic_Lyapunov_Functions}. \QEDB 
\end{remark}

The previous facts and discussion motivate the main problem statement considered in this section:

\noindent 
\textbf{Problem 1:} Under the Assumption that the HDS \eqref{eq:combined_org_HDS} is strongly $\nabla V$-stabilizable, design an implementable feedback law $u$ that renders UGAS the compact set $\mathcal{A}$.

\subsection{Robust Model-Free Stabilization via Oscillatory Hybrid Feedback Control}
To solve Problem 1, we propose a class of hybrid \emph{model-free} feedback controllers that require only real-time measurement or evaluations of the SCLF $V(x)$. 
Specifically, for $i\in\{1,2,\dots,r\}$, we propose the following control law:
\begin{subequations}\label{eq:ES_Law}
\begin{align}
    u^\varepsilon_i(V(x),\eta) = \varepsilon^{-1}\sqrt{\frac{4\pi\gamma}{T_i\kappa}} \langle\exp(\kappa V(x) S)e_1,\eta_i \rangle,
\end{align}
where $V$ is a SCLF with respect to $\mathcal{A}$ for $\mathcal{H}$, the constants $\kappa,\gamma,\varepsilon \in\mathbb{R}_{>0}$ are tuning parameters, $\eta_i\in\mathbb{S}^1$ for $i\in\{1,2,\dots,r\}$ is the state of the linear oscillator
\begin{align}\label{eq:single_oscillator_dynamics}
\eta_i\in\mathbb{S}^1,~~\dot{\eta}_i&= 2\pi T_i^{-1}\varepsilon^{-2} S \eta_i, & S&= \begin{pmatrix}
        0 & 1\\-1 & 0
    \end{pmatrix},
\end{align}
where $\mathbb{S}^1=\{\eta_i\in\mathbb{R}^2:\eta_{i,1}^2+\eta_{i,2}^2=1\}$ is the unitary circle,
and $\{T_1,T_2,\dots,T_r\}\subset \mathbb{Q}_{>0}$ is a collection of constants satisfying $T_i \neq T_j$ for all $i\neq j$. For brevity, we will use
\begin{align}\label{eq:lambda}
    \dot{\eta}&= \Lambda_\varepsilon(\eta),
\end{align}
\end{subequations}
to denote the collective (continuous-time) dynamics of the uncoupled oscillators $\eta=(\eta_1,\eta_2,\dots,\eta_r)\in\mathbb{S}^1\times\mathbb{S}^1\times\cdots\times\mathbb{S}^1=:\mathbb{T}^r$, where the map $\Lambda_\varepsilon$ on the right hand side of \eqref{eq:lambda} is defined consistently with \eqref{eq:single_oscillator_dynamics}. 

By applying the proposed feedback law \eqref{eq:ES_Law} to close the loop for $\mathcal{H}$ in \eqref{eq:combined_org_HDS}, the resulting closed-loop HDS has the state $y = (\xi, \eta)$ and dynamics
\begin{subequations}\label{eq:closed_loop_org_HDS}
\begin{align}
    \mathcal{H}_{cl}:\begin{cases}
            ~y\in C\times\mathbb{T}^r,&
            \dot{y}\hphantom{^+} \in \hat{F}_\varepsilon(y):=F_\varepsilon(\xi,\eta)\times F_e(\theta)
        \\
        ~y\in D\times\mathbb{T}^r, & y^+\in \hat{G}(y):=G(\xi)\times\{\eta\}
        \end{cases}
\end{align}

\vspace{-0.5cm}
\noindent 
where, for each $\varepsilon\in\mathbb{R}_{>0}$, the set-valued map $F_\varepsilon$ is given by
\begin{align}
    F_\varepsilon(\xi,\eta)&= \{f_\varepsilon(\xi,\eta))\}\times \{\Lambda_\varepsilon(\eta)\}\times F_e(\theta),\label{fsubespilon}\\
    f_\varepsilon(\xi,\eta)&= f_0(x,\theta) + \sum_{i=1}^r f_i(x,\theta) u^\varepsilon_i(V(x),\eta).
\end{align}
\end{subequations}
Using $\bar{\mathcal{A}}:=\mathcal{A}\times\Theta\times \mathbb{T}^r$, we can now state our first main result. The proof is presented in Section \ref{sec:Proofs}.
\begin{thm}\normalfont\label{thm:main_theorem}
    Let $\mathcal{H}$ in \eqref{eq:combined_org_HDS} be strongly $\nabla V$-stabilizable. Then, $\bar{\mathcal{A}}$ is SGpAS as $\varepsilon\rightarrow 0^+$ for the closed-loop  $\mathcal{H}_{cl}$ in \eqref{eq:closed_loop_org_HDS}. 
\end{thm}
We also have the following corollary, which is of independent interest for problems defined on compact spaces.
\begin{cor}\label{cor:main_theorem_corollary}
     Let $\mathcal{H}$ in \eqref{eq:combined_org_HDS} be strongly $\nabla V$-stabilizable. If $\mathcal{X}$ is compact, then $\bar{\mathcal{A}}$ is UGpAS as $\varepsilon\rightarrow 0^+$ for the closed-loop HDS $\mathcal{H}_{cl}$ in \eqref{eq:closed_loop_org_HDS}.
\end{cor}

The proof of Theorem \ref{thm:main_theorem}, presented in Section \ref{sec:Proofs}, exploits the Lie-Bracket averaging theorems for well-posed HDS introduced in \cite{abdelgalil2023lie}. However, in contrast to the results in \cite{abdelgalil2023lie}, which consider time-varying HDS, in the present paper we modeled the complete dynamics as time-invariant HDS by using the dynamic oscillators \eqref{eq:lambda}, which evolve in the compact set $\mathbb{T}^r$. Similar oscillators have been considered before for the study of hybrid extremum-seeking control via first-order averaging \cite{PoTe17Auto,PovedaNaliAuto20}. However, Theorem \ref{thm:main_theorem} is the first result that uses time-invariant oscillators for the analysis of Lie-bracket-based averaging algorithms.
\begin{remark}\normalfont
The results of Theorem \ref{thm:main_theorem} extend the model-free control laws studied for ODEs in \cite{scheinker2017model}, \cite[Sec. 6]{scheinker2024100} and, more recently, in \cite{wang2024extremum}, to the framework of hybrid systems.
\end{remark}
The following important corollary is a consequence of \cite[Thm. 7.21]{bookHDS} and $\mathcal{H}_{cl}$ satisfying the Hybrid Basic Conditions of Assumption \ref{asmp:basic}.
\begin{cor}[Robustness]\normalfont
    Consider the perturbed closed-loop HDS $\mathcal{H}^d_{cl}$ obtained from \eqref{eq:closed_loop_org_HDS} as follows: 
    \begin{align*}
    \mathcal{H}^{d^*}_{cl}:\begin{cases}
            ~~y+d_1\in C\times\mathbb{T}^r,&
            \dot{y}\hphantom{^+} \in \hat{F}_\varepsilon(y+d_2)+d_3,
        \\
        ~~y+d_4\in D\times\mathbb{T}^r, & y^+\in \hat{G}(y+d_5)+d_6,
        \end{cases}
    \end{align*}
    where the signals $d_i:\text{dom}(y)\rightarrow (C\cup D)\times\mathbb{T}^r$ are measurable functions satisfying $\sup_{(t,j)\in\text{dom}(y)}|d_i(t,j)|\leq d^*$,
    for all $i\in\{1,2,\dots,6\}$ and some $d^*\in\mathbb{R}_{>0}$. Then, the perturbed HDS $\mathcal{H}^d_{cl}$ renders the set $\bar{\mathcal{A}}$ SGpAS as $(d^*,\varepsilon)\rightarrow 0^+$. Moreover, if $\mathcal{X}$ is compact, then the set $\bar{\mathcal{A}}$ is UGpAS as $(d^*,\varepsilon)\rightarrow 0^+$ for system $\mathcal{H}^d_{cl}$.
\end{cor}
In the next section, we present different applications of Theorem \ref{thm:main_theorem} in the context of robust global stabilization of a point $p^\star$ on a smooth manifold $\mathcal{M}$ that is not globally contractible, and under unknown switching control directions.

\section{Minimum Seeking for Synergistic Potential Functions on Smooth Manifolds}\label{sec:Synergistic_Lyapunov_Functions}
Let $\mathcal{M}$ be a smooth closed manifold properly embedded within, and equipped with the Riemannian metric of, an ambient Euclidean space $\mathbb{R}^{n_p}$. Let $p\in\mathcal{M}$ denote the state of a plant with dynamics
\begin{align}\label{eq:manifold_control_affine}
    \dot{p}&= \sum_{i=1}^r b_i(p) \theta_i u_i,
\end{align}
where $u=(u_1,u_2,\dots,u_r)$ are the control inputs. The vector $(\theta_1,\theta_2,\dots,\theta_r)$ corresponds to unknown control gains that are allowed to dynamically switch between finitely many values, i.e. $(\theta_1,\theta_2,\dots,\theta_r)\in \mathcal{E}\subset\mathbb{R}^r$ and $|\mathcal{E}|< +\infty$. For example, without loss of generality, we may consider the normalized directions $\mathcal{E}=\{+1,0,-1\}^r$. 
\begin{remark}\normalfont
For the purpose of illustration, in \eqref{eq:manifold_control_affine} we have specialized the form of uncertainty in the control directions to be linear. Nevertheless, we emphasize that general, potentially nonlinear, dependence on $\theta$ is allowed as long as Assumptions \ref{asmp:regularity_org_HDS} and \ref{asmp:exogeneous_signal_regularity} hold. 
\end{remark}

To satisfy Assumption \ref{asmp:exogeneous_signal_regularity}, we impose a \emph{dwell time condition} on the rate of switching of $\theta$. We also require that each control gain $\theta_i$ does not vanish for an indefinite duration of time, since otherwise uniform (practical) asymptotic stability properties will be precluded. The two requirements are satisfied if the vector $(\theta_1,\theta_2,\dots,\theta_r)$ is governed by the HDS \eqref{eq:exogenous_signal_dynamics} with state $\theta=(\theta_1,\theta_2,\dots,\theta_r,\theta_{r+1},\theta_{r+2})\in\Theta_C\cup \Theta_D$, and the following data:
\begin{subequations}\label{hybridautomaton}
    \begin{align}
        \Theta_C&:= \mathcal{E}\times [0,1]\times[0,T_\circ],\\
        \Theta_D&:= \mathcal{E}\times \{1\}\times[0,T_\circ], \\
        F_e(\theta)&:= \{0\}\times \left[0,\chi_1\right]\times(\left[0,\chi_2\right]-\mathbb{I}_{\mathcal{E}_b}(\theta)), \\
        \mathcal{E}_b&:=\{\theta ~|~ \exists i\in\{1,2,\dots,r\} \text{ s.t. } \theta_i = 0\},\\
      \tilde{G}_e(\theta)&:= \mathcal{E}\backslash\{(\theta_1,\theta_2,\dots,\theta_r)\},\\
        G_e(\theta)&:= \tilde{G}_e(\theta) \times \{0\}\times\{\theta_{r+2}\},
    \end{align}
\end{subequations}
for some $T_\circ\in\mathbb{R}_{>0}$, $\chi_1\in\mathbb{R}_{>0}$, $\chi_2\in(0,1)$, and where $\mathbb{I}_{\mathcal{E}_b}(\cdot)$ is the classical indicator function on the subset $\mathcal{E}_b$. 

The construction \eqref{hybridautomaton} involves a hybrid automaton, similar to the one considered in \cite[Prop. 1.1]{ConverseHybrid2}, with auxiliary state $\theta_{r+1}\in\mathbb{R}_{\geq0}$, and a time-ratio monitor with auxiliary state $\theta_{r+2}\in\mathbb{R}_{\geq0}$, similar to the one considered in \cite[Lemma 7]{PoTe17Auto}. In particular, every time the condition $\theta_{r+1}=1$ is satisfied, $\theta_{r+1}$ is reset to zero, and the vector of directions $\theta$ is allowed to jump according to the rule $\theta^+\in \tilde{G}_e(\theta)$, while $\theta_{r+2}$ remains unchanged. The role of $\theta_{r+2}$ is to model a monitor that decreases continuously whenever the current direction $\theta$ has a component equal to zero, i.e., when $\theta\in\mathcal{E}_b$ and therefore at least one of the entries in \eqref{eq:manifold_control_affine} has zero controllability. In this case, the condition $\theta_{r+2}=0$ will eventually occur, forcing $\theta$ to jump out of the set $\mathcal{E}_b$ and into the set $\{-1,1\}^r$. Note that, unlike existing hybrid monitors for switching systems with unstable modes \cite{PoTe17Auto}, the set $\{-1,1\}^r$ also contains vectors $\theta$ with negative signs, which could be destabilizing for system \eqref{eq:manifold_control_affine} under a feedback law designed for $\theta=\mathbf{1}_r$. In other words, in this paper, the time-ratio monitor prevents the system from spending too much time in modes with zero control directions but not in modes with potentially problematic negative control directions.

By construction, and by \cite[Prop. 1.1]{ConverseHybrid2}, and \cite[Lemma 7]{PoTe17Auto}, each hybrid arc generated by the HDS defined by \eqref{eq:exogenous_signal_dynamics} with data \eqref{hybridautomaton} satisfies the following dwell-time and activation-time inequalities for any two times $t_2>t_1$ in its domain:
\begin{subequations}
\begin{align}
    N_\sharp(t_1,t_2)&\leq \chi_1(t_2-t_1) + 1, \label{ADT}\\
    T_\sharp(t_1,t_2)&\leq \chi_2(t_2-t_1) + T_\circ,\label{ATT}
\end{align}
\end{subequations}
where $N_\sharp(t_1,t_2)$ is the total number of jumps during the time interval $(t_1,t_2)$, and
\begin{align}
    T_\sharp(t_1,t_2)=\int_{t_1}^{t_2}\mathbb{I}_{\mathcal{E}_b}(\theta(t,j))\,\text{d}t,
\end{align}
is the total time that the unknown vector $\theta$ spends in the set $\mathcal{E}_b$ during the same interval \cite{PoTe17Auto}.

Next, we impose the following standard regularity condition on the control vector fields $b_i$ of system \eqref{eq:manifold_control_affine}. This condition essentially guarantees that, in the ideal setting when $\theta=\mathbf{1}_r$, system \eqref{eq:manifold_control_affine} can actually be controlled by having vector fields $b_i$ that span the tangent space of the manifold $\mathcal{M}$ everywhere. Such assumptions are standard in the literature related to model-free stabilization and optimization \cite{DurrManifold,AlexTAC}. Similar assumptions have also been used in the design of synergistic hybrid feedback laws that overcome topological obstructions to stability \cite{TrackingManifold}. 
\begin{asmp}\label{asmp:synergistic_manifold_regularity}
    For all $i\in\{1,2,\dots,r\}$ and all $p\in\mathcal{M}$, $b_i$ is $\mathcal{C}^1$, $b_i(p)\in T_p\mathcal{M}$, and there exists a constant $\lambda>0$ such that for all all $v\in T_p\mathcal{M}$, it holds that $\textstyle\sum_{i=1}^r \langle b_i(p), v\rangle^2 \geq \lambda \langle v, v\rangle$. \QEDB 
\end{asmp}
%


Let $V$ be a $\mathcal{C}^1$ function, and consider the subset
\begin{align}
    \text{Crit}(V):=\{p\in\mathcal{M}~|~\nabla V(p)\in (T_p\mathcal{M})^\perp\},
\end{align}
where $(T_p\mathcal{M})^\perp\subset\mathbb{R}^{n_p}$ denotes the orthogonal complement of the tangent space $T_p\mathcal{M}$ with respect to the Riemmannian metric of the ambient Euclidean space $\mathbb{R}^{n_p}$. With this notation, we recall the following definition, adapted from \cite{Mayhew10Thesis}.
\begin{definition}\normalfont \label{defn:syn_lyap_func}
    For $N\in\mathbb{N}_{\geq 1}$ and $\mathcal{Q}=\{1,2,\dots,N\}$, a family of $\mathcal{C}^1$ functions $\{V_q\}_{q\in\mathcal{Q}}$, $V_q:\mathcal{M}\rightarrow\mathbb{R}_{\geq 0}$, is said to be a \emph{$\delta$-gap synergistic family of potential functions with respect to $p^\star$} if, for all $q\in\mathcal{Q}$, the following holds:
    \begin{enumerate}
        \item $V_q$ is positive definite with respect to $p^\star$;
        \item $\forall c\in\mathbb{R}_{\geq 0}$, $\{ p\in\mathcal{M}~|~V_q(p)\leq c\}$ is compact;
        \item There exists an open neighborhood $\mathcal{U}_q\subset\mathbb{R}^{n_p}$ of $p^\star$ such that $\mathcal{U}_q\cap \text{Crit}(V_q)=p^\star$;
        \item There exists $\delta\in\mathbb{R}_{>0}$ such that $\delta <\Delta^\star$, where $\Delta^\star$ is
        \begin{align}
            \Delta^\star=\min_{\substack{q\in\mathcal{Q}\\ p\in \text{Crit}(V_q)}} V_q(p) - \max_{\tilde{q}\in\mathcal{Q}} V_{\tilde{q}}(p).
        \end{align}
    \end{enumerate}
\end{definition}
The construction of synergistic families of potential functions for stability problems typically requires qualitative information on the underlying manifold $\mathcal{M}$ and the target point $p^\star \in\mathcal{M}$. We make the assumption that such functions are available to us, and later we show how to satisfy this assumption in different applications.
\begin{asmp}\label{asmp:synergistic_family}
The family of functions $\{V_q\}_{q\in\mathcal{Q}}$, $\mathcal{Q}=\{1,2,\dots,N\}$, is a $\delta$-gap synergistic family of potential functions with respect to $p^\star$. 
\end{asmp}

To exploit the existence of synergistic potential functions, we introduce a logic state $q\in\mathcal{Q}$ into our controller, and we let $x=(p,q)\in\mathcal{M}\times\mathcal{Q}$, and define
\begin{subequations}\label{eq:synergistic_mnfd_data}
\begin{align}
    \mathcal{X}_C&:=\{x\in\mathcal{M}\times\mathcal{Q}~|~\mu(x)\leq \delta \},\\
    \mathcal{X}_D&:=\{x\in\mathcal{M}\times\mathcal{Q}~|~\mu(x)\geq \delta \},\\
    G_{\mathcal{X}}(x)&:= \{p\}\times
        \big\{\arg\min_{\tilde{q}\in\mathcal{Q}}V(p,\tilde{q})\big\},
\end{align}
where the function $\mu:\mathcal{M}\times\mathcal{Q}\rightarrow \mathbb{R}_{\geq 0}$ is defined by
\begin{align}
    \mu(x):= V_q(p) - \textstyle\min_{\tilde{q}\in\mathcal{Q}} V_{\tilde{q}}(p).
\end{align}
Also, let $f_0$ and $f_i$ for $i\in\{1,2,\dots,r\}$ be given by
\begin{align}
    f_0(x,\theta) &= 0, & f_i(x,\theta)&:=\theta_i\begin{pmatrix} b_i(p) \\ 0 \end{pmatrix}.
\end{align}
\end{subequations}
where the second entry in $f_i$ models the continuous-time dynamics of the logic state $q$, which remains constant during flows, i.e., $\dot{q}=0$. Finally, define the function $V:\mathcal{M}\times\mathcal{Q}\rightarrow\mathbb{R}_{\geq 0}$ as follows:
\begin{equation}\label{importantCLF}
V(x) = V_q(p).
\end{equation}
Then, we have the following proposition that applies to any system of the form \eqref{eq:manifold_control_affine} that satisfies the previous assumptions. The proof is presented in Section \ref{sec:Proofs}.
\begin{prop}\normalfont \label{prop:SCLF_synergistic}
      Suppose that Assumptions \ref{asmp:synergistic_manifold_regularity} and \ref{asmp:synergistic_family} are satisfied, and $V$ is given by \eqref{importantCLF}. Then, $\mathcal{H}$ is strongly $\nabla V$-stabilizable. \QEDB 
\end{prop}
As a consequence of Theorem \ref{thm:main_theorem} and Proposition \ref{prop:SCLF_synergistic}, we are able to conclude that, under the feedback law \eqref{eq:ES_Law}, the subset $\mathcal{A}\times\Theta\times\mathbb{T}^r$ is SGpAS as $\varepsilon\rightarrow 0^+$ for the closed loop HDS $\mathcal{H}_{cl}$ \eqref{eq:closed_loop_org_HDS}.

\begin{figure}[t]
    \centering
\includegraphics[width=\linewidth]{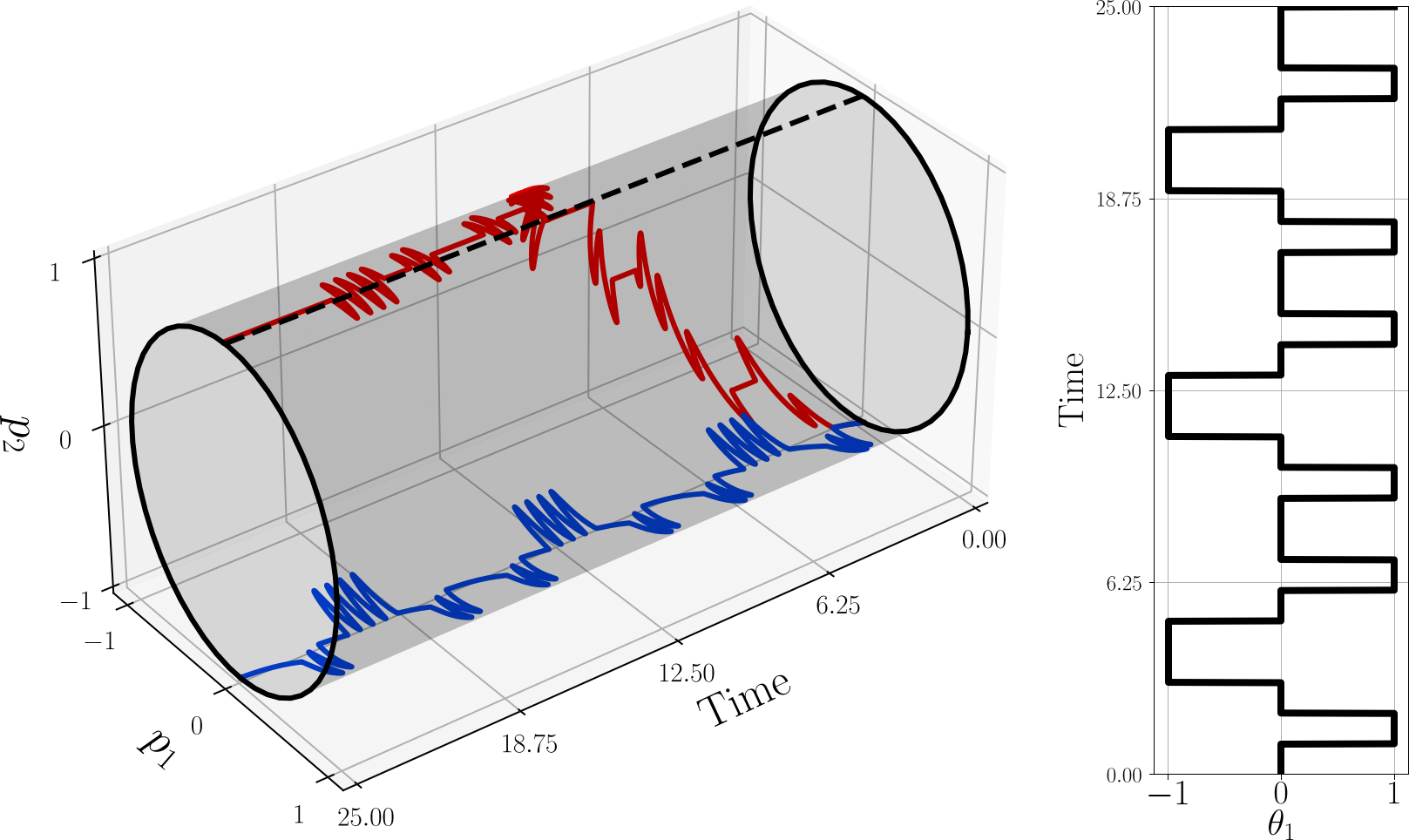}
    \caption{Global stabilization of $p^\star=(0,1)$ on $\mathbb{S}^1$ under unknown  switching control directions. Left: trajectories of \eqref{systemonthecircle} under the proposed model-free feedback law (red color), and the non-hybrid feedback law \cite{scheinker2017model} (blue color). Right: Evolution in time of $\theta_1$.}
    \label{fig:circle}
\end{figure}

Next, we present several novel concrete applications of the results in this section along with numerical simulation results. Henceforth, we take $\mathcal{E}=\{+1,0,-1\}^r$. Since the construction of synergistic families of potential functions is not the main contribution of our manuscript, we rely on existing results in the literature \cite{Mayhew10Thesis}. Instead, we focus on the role of the adaptive feedback law \eqref{eq:ES_Law} and its ability to ``emulate'' the ideal control law \eqref{eq:avg_feedback_law} in the context of hybrid control using highly oscillatory feedback.
\subsection{Robust Global Stabilization on $\mathbb{S}^1$}
As the first example, we take $\mathcal{M}=\mathbb{S}^1$, $r=1$, and we consider the control-affine system
\begin{align}\label{systemonthecircle}
    \dot{p}&= b_1(p)\theta_1 u_1, & b_1(p)&= S p,
\end{align}
where $p\in\mathbb{S}^1$, $S$ is the matrix defined in \eqref{eq:ES_Law}, $\theta_1\in\mathcal{E}=\{+1,0,-1\}$ is the unknown control gain, and $u_1$ is the control input.  The goal is to \emph{globally} stabilize a point $\vartheta^\star\in\mathbb{S}^1$. This problem was solved in \cite{mayhew2010hybrid} and \cite{ochoa2022model} under the assumption of having constant and known control directions $\theta_1:=1$. 

To globally stabilize $\vartheta^\star\in\mathbb{S}^1$ under unknown control directions, we consider the synergistic family of potential functions $\{W_1,W_2\}$ defined by
\begin{subequations}\label{eq:syn_family_circle}
    \begin{align}
        W_q(p) &:= W\circ \Phi_q(p), & W(p)&:= 1-\langle \vartheta^\star,p\rangle,
    \end{align}
    where $q\in\mathcal{Q}=\{1,2\}$, and the maps $\Phi_q:\mathbb{S}^1\rightarrow\mathbb{S}^1$ are
    \begin{align}
        \Phi_q(p):=\exp\left((3/2-q) W(p) S\right)p.
    \end{align}
\end{subequations}
As shown in \cite{mayhew2010hybrid,Mayhew10Thesis}, the family of functions $\{W_1,W_2\}$ is a $\delta$-gap synergistic family of potential functions function with respect to $\vartheta^\star$ for any $\delta\in(0,1)$. It follows from Proposition \ref{prop:SCLF_synergistic} that the function $V$ defined by $V(p,q)=W_q(p)$ is an SCLF with respect to $\mathcal{A}=\{p^\star\}\times\mathcal{Q}$ for $\mathcal{H}$, and therefore that $\mathcal{A}\times\Theta\times\mathbb{S}^1$ is SGpAS for $\mathcal{H}_{cl}$. Due to the compactness of $\mathbb{S}^1$, we can invoke Corollary \ref{cor:main_theorem_corollary} to conclude that $\mathcal{A}\times\Theta\times\mathbb{S}^1$ is UGpAS.

Numerical simulations results illustrating the performance of the proposed controller for this example are shown in Figure \ref{fig:circle}. To generate the results, we used $\gamma = \sqrt{1}$, $\kappa = 4$, $\varepsilon=1/\sqrt{4\pi}\approx 0.28$. We also used $\delta = 1/4$ for the synergistic family of potential functions $\{W_1,W_2\}$. The target point is $p^\star=(0,1)$. The right plot shows the evolution in time of the control direction $\theta$, which vanishes during bounded (but persistent) periods of time. Finally, to emphasize the robustness of the proposed feedback law, we added a small adversarial perturbation that locally stabilizes (in the absence of switching) the problematic critical point $p^\sharp=(0,-1)$. As shown in the figure, the proposed control law is not affected by the perturbation whereas a non-hybrid model-free feedback law is effectively trapped by the adversarial perturbation in the vicinity of the critical point $p^\sharp$.
\begin{figure}[t]
    \centering
    \includegraphics[width=\linewidth]{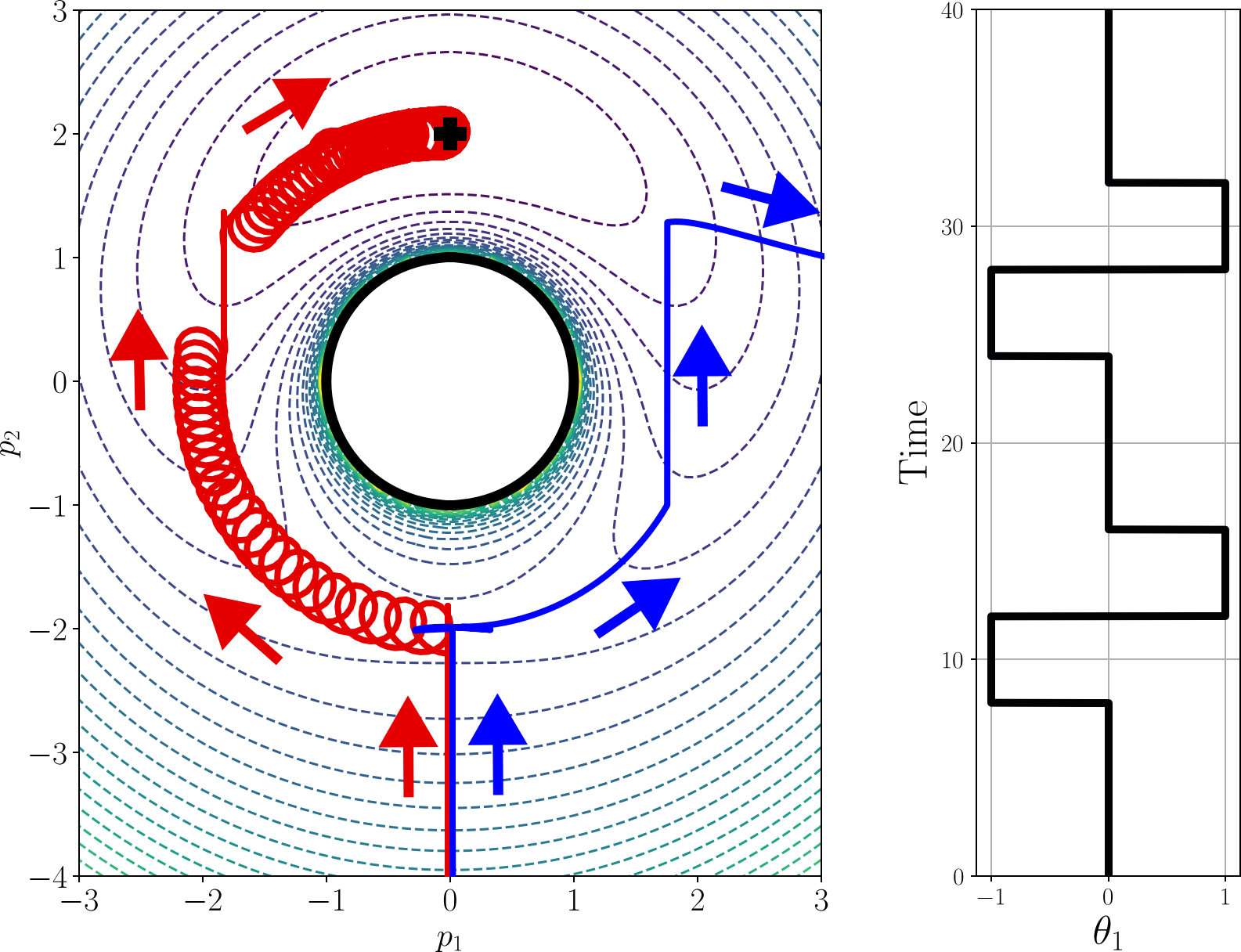}
    \caption{Target-seeking with obstacle avoidance under unknown control gains. Left: trajectories of the vehicle under the proposed hybrid model-free feedback law (shown in red), and under the vanilla synergistic hybrid feedback (shown in blue) \cite{TrackingManifold}. Right: Evolution in time of the control gain $\theta_1$. The gain $\theta_2=1$ is taken to be constant.}
    \label{fig:obstacle_avoidance}
\end{figure}
\subsection{Robust Target Seeking with Obstacle Avoidance}\label{subsec:obstacle_avoidance_integrator}
Consider the problem of stabilizing a target position for a mobile vehicle moving in an obstructed planar domain. Let $\mathcal{O}\subset\mathbb{R}^2$ denote the obstacle and let $z=(z_1,z_2)\in\mathbb{R}^2\backslash\mathcal{O}$ denote the position of the vehicle. The goal of the vehicle is to reach a target position $z^\star=(z_1^\star,z_2^\star)\in\mathbb{R}^2\backslash\mathcal{O}$ while avoiding the obstacle $\mathcal{O}$. We assume that the motion of the vehicle is governed by the kinematic equations
\begin{align}
    \dot{z}&= \sum_{i=1}^2 e_i\,\theta_i \,u_i, 
\end{align}
where $u=(u_1,u_2)\in\mathbb{R}^2$ are the control input, and $(\theta_1,\theta_2)\in \mathcal{E}=\{+1,0,-1\}^r$ are the unknown control gains. Since the obstacle $\mathcal{O}$ is bounded, there exists $z_{\mathcal{O}}\in\mathbb{R}^2$ and $d_0 > 0$ such that $\mathcal{O}\subset z_{\mathcal{O}}+d\,\mathbb{B}$ for all $d\geq d_0$. To guarantee feasibility, we impose the following assumption on the target point $z^\star$.
\begin{asmp}\label{asmp:obstacle_contained}
    $\exists~d^\star > d_0$ such that $z^\star\in \mathbb{R}^2\backslash(z_{\mathcal{O}}+d^\star\,\mathbb{B})$.
\end{asmp}
To avoid the obstacle, and following the ideas of \cite{TrackingManifold}, we consider the map $\varphi:\mathbb{R}^2\backslash(z_{\mathcal{O}}+d_0\,\mathbb{B}))\rightarrow \mathbb{R}\times\mathbb{S}^1$ given by 
\begin{align}\label{eq:obstacle_avoidance_diffeomorphism}
    \varphi(z) &= (\log(|z-z_{\mathcal{O}}|-d^\star),(z-z_{\mathcal{O}})/|z-z_{\mathcal{O}}|),
\end{align}
As shown in \cite{TrackingManifold}, $\varphi$ is a well-defined diffeomorphism. The \emph{pushforward} \cite{lee2012smooth} of the kinematics of the vehicle under the diffeomorphism $\varphi$ are given by
\begin{align}
    \dot{p}&= \sum_{i=1}^2 b_i(p)\theta_i u_i, & b_i(p)&=\text{D}\varphi\circ\varphi^{-1}(p) e_i,
\end{align}
where $p=(\rho,\vartheta)\in\mathbb{R}\times\mathbb{S}^1$. Therefore, global stabilization of the target position $z^\star$ in $\mathbb{R}^2\backslash(p_{\mathcal{O}}+d_0\,\mathbb{B}))$ is equivalent to globally stabilizing the point $p^\star=(\rho^\star,\vartheta^\star)$ on the smooth manifold $\mathcal{M}=\mathbb{R}\times\mathbb{S}^1$. However, the topology of $\mathcal{M}$ prohibits global stabilization by continuous feedback since $\mathbb{S}^1$ is a compact boundary-less manifold. With that in mind, we introduce the family of functions
\begin{align*}
    V_q(p)&:= \frac{1}{2}(\rho-\rho^\star)^2 + \sqrt{(\text{e}^{\rho}-\text{e}^{\rho^\star})^2+1}-1 + W_q(\vartheta),
\end{align*}
where $q\in\mathcal{Q}=\{1,2\}$, and $W_q$ are the functions defined in \eqref{eq:syn_family_circle}. As in the previous subsection, the family of functions $\{W_q\}_{q\in\mathcal{Q}}$ is a $\delta$-gap synergistic Lyapunov function for $\vartheta^\star$ on $\mathbb{S}^1$ for any $\delta\in(0,1)$, which follows by \cite{Mayhew10Thesis}. Therefore, it is straightforward to show that the family of functions $\{V_1,V_2\}$ is a $\delta$-gap synergistic family of potential functions with respect to $p^\star$. Thus, from Proposition \ref{prop:SCLF_synergistic}, the function $V$ defined by $V(p,q)=V_q(p)$ is an SCLF with respect to $\mathcal{A}=\{p^\star\}\times \mathcal{Q}$ for $\mathcal{H}$. By invoking Theorem \ref{thm:main_theorem}, we conclude that $\mathcal{A}\times\Theta\times\mathbb{T}^2$ is SGpAS for the original hybrid system with unknown and dynamic control directions $\mathcal{H}_{cl}$. 

To demonstrate the performance of the proposed controller compared to existing synergistic hybrid feedback controllers \cite{TrackingManifold}, we present numerical simulations in Figure \ref{fig:obstacle_avoidance}. To generate the results, we used $\gamma = 2$, $\delta = 1/4$, $\kappa = 4$, and $\varepsilon=1/\sqrt{6\pi}\approx 0.165$. The target position is $z^\star=(0,2)$, and an obstacle with radius $d=1$ is centered at the origin, i.e. $z_{\mathcal{O}}=(0,0)$.  Due to the special structure of the control vector fields in this example, we are able to use a single oscillator $\eta_1$ with period $T_1=1$ and rely on a $\pi/2$ phase shift to guarantee non-resonance between the inputs $u_1$ and $u_2$. Figure \ref{fig:obstacle_avoidance} clearly indicates that the vanilla synergistic hybrid feedback fails to reach the target. This is to be expected since traditional synergistic hybrid controllers require knowledge of the control direction. By contrast, the proposed algorithm is able to reach the target despite the persistent switching in the control gain sign and the temporary loss of control.
\begin{figure}[t]
    \centering
    \includegraphics[width=\linewidth]{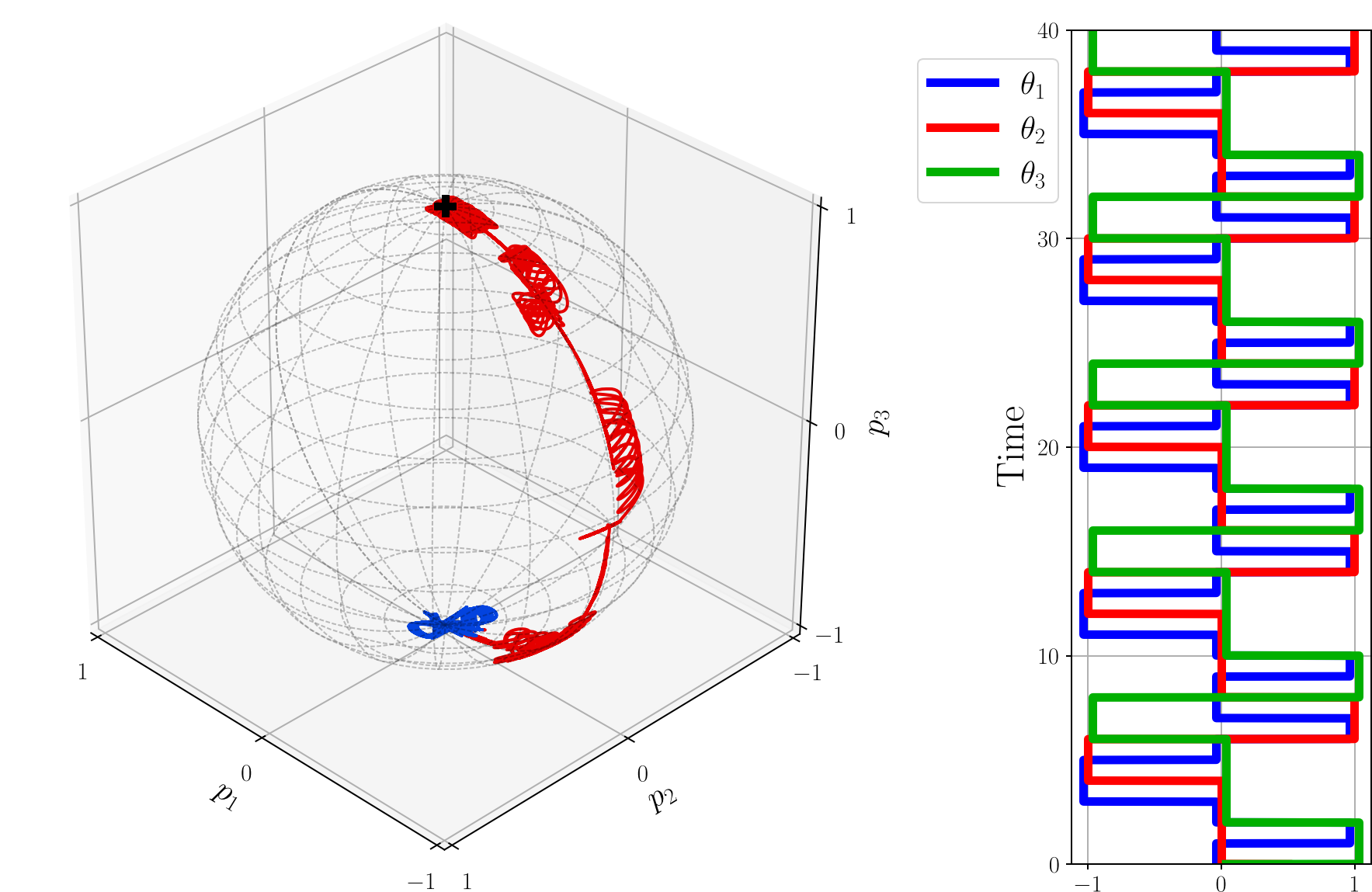}
    \caption{Global stabilization of $p^\star=(0,0,1)$ on $\mathbb{S}^2$ under unknown persistently switching control directions. Left: trajectory corresponding to the proposed hybrid model-free feedback law \eqref{eq:ES_Law} (shown in red), and the non-hybrid, model-free feedback law \cite{scheinker2017model,durrEbenbauer} (shown in blue). Right: Evolution in time of the control gains $(\theta_1,\theta_2,\theta_3)$.}
    \label{fig:sphere}
\end{figure}
\subsection{Robust Global Stabilization on $\mathbb{S}^2$}

\vspace{-0.3cm}\noindent
For the third example, we take $\mathcal{M}=\mathbb{S}^2$, $r=3$, and we consider the control-affine system

\vspace{-0.4cm}\noindent
\begin{align}
    \dot{p}&= \sum_{i=1}^3 b_i(p)\theta_i u_i, & b_i(p)&= e_i-\langle p,e_i\rangle p,
\end{align}

\vspace{-0.4cm}\noindent
where $p\in\mathbb{S}^2$, $(\theta_1,\theta_2,\theta_3)\in \mathcal{E}=\{+1,0,-1\}^3$ are the unknown control gains, and $u=(u_1,u_2,u_3)\in\mathbb{R}^3$ are the control inputs. The goal is to globally stabilize an arbitrary point $p^\star\in\mathbb{S}^2$. To that end, we introduce the synergistic family of potential functions
\begin{align}
    V_q(p)&:= W\circ \Phi_q(p), & W(p)&:= 1-\langle p, p^\star\rangle,
\end{align}
where $q\in\{1,2\}$, the maps $\Phi_q:\mathbb{S}^2\rightarrow\mathbb{S}^2$ are defined by
\begin{align}
    \Phi_q(p):=\text{exp}\left((3/2-q) W(p) [p^\star_\perp]_{\times}\right)p,
\end{align}
\begin{figure*}[t!]
    \centering
    \includegraphics[width=\linewidth]{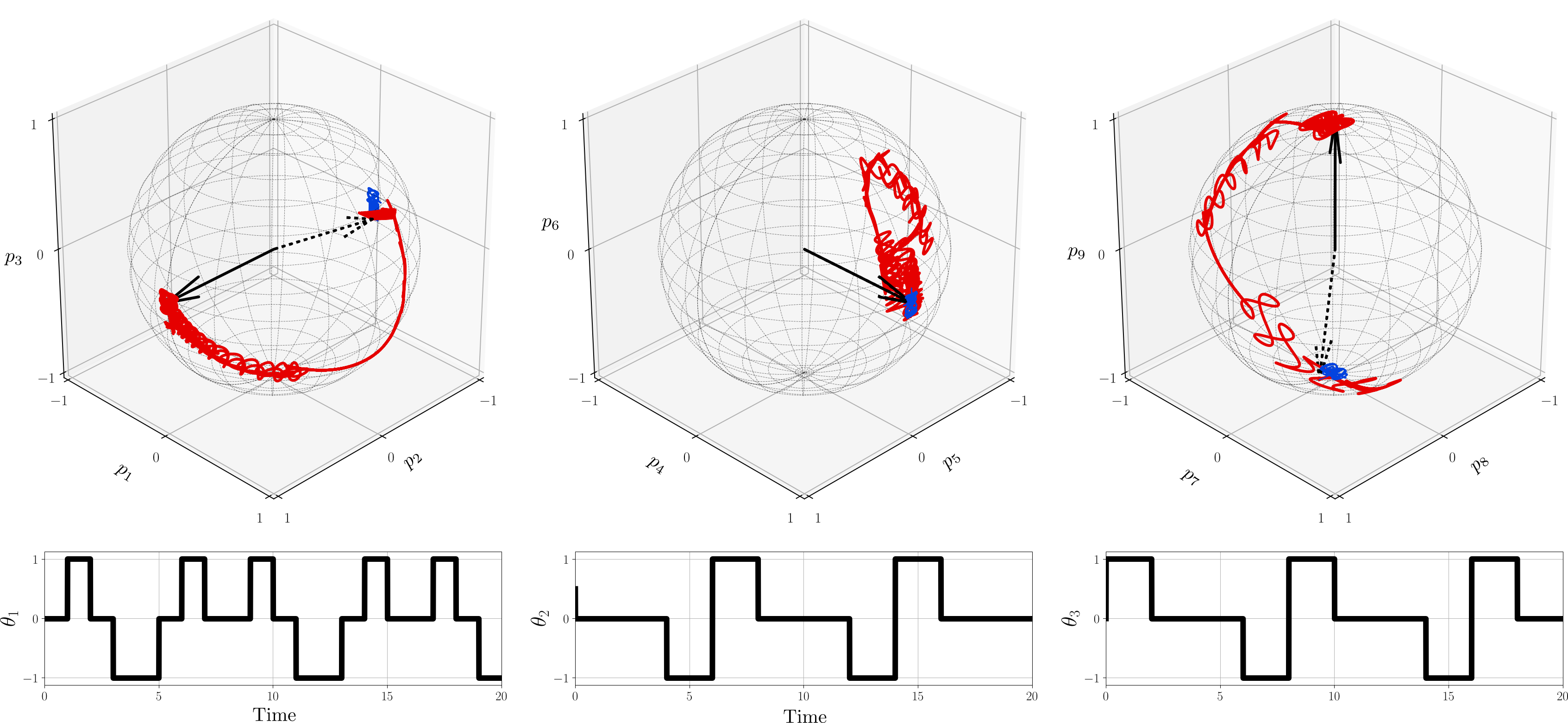}
    \caption{Global stabilization of $p^\star=vec(I)$ on $\text{SO}(3)$ under unknown switching control directions (shown in the bottom plot). The red trajectory corresponds to the proposed model-free feedback law \eqref{eq:ES_Law}, whereas the blue trajectory corresponds to a non-hybrid model-free feedback law \cite{scheinker2017model}. Both systems are subject to the same small bounded persistent disturbance, which disrupts the blue trajectory.}
    \label{fig:rigid_body}
\end{figure*}
and $p^\star_\perp\in\mathbb{S}^2$ is such that $\langle p^\star_\perp,p^\star\rangle = 0$. As shown in \cite{Mayhew10Thesis}, the family of functions $\{V_1,V_2\}$ is a $\delta$-gap synergistic family of potential functions function with respect to $p^\star$ for any $\delta\in(0,1)$. It follows from Proposition \ref{prop:SCLF_synergistic} that the function $V$ defined by $V(p,q)=V_q(p)$ is an SCLF with respect to $\mathcal{A}=\{p^\star\}\times\mathcal{Q}$ for $\mathcal{H}$, and therefore that $\mathcal{A}\times\Theta\times\mathbb{S}^2$ is SGpAS for $\mathcal{H}_{cl}$. By invoking Corollary \ref{cor:main_theorem_corollary}, we conclude that $\mathcal{A}\times\Theta\times\mathbb{T}^3$ is UGpAS.

\vspace{-0.1cm}\noindent
Numerical simulations results are shown in Figure \ref{fig:sphere}. To generate the results, we used $\gamma = \sqrt{1}$, $\kappa = 4$, $T_1=3$, $T_2=2$, $T_3=1$, $\delta = 1/5$, and $\varepsilon=1/\sqrt{8\pi}$. The target point is $p^\star=(0,0,1)$, and $p_\perp^\star=(0,1,0)$. Finally, to emphasize the robustness of the proposed feedback law compared to the standard non-hybrid model-free controllers \cite{DurrManifold}, we added a small adversarial perturbation that locally stabilizes the problematic critical point $p^\sharp=(0,0,-1)$ in the absence of switching. As shown in Figure \ref{fig:sphere}, the proposed control law is not affected by the perturbation whereas a vanilla model-free feedback law is trapped by the adversarial perturbation in the vicinity of the critical point $p^\sharp$.
\subsection{Robust Global Stabilization on $\text{SO}(3)$}

\vspace{-0.2cm}\noindent
Next, we consider the problem of globally stabilizing a desired attitude $R^\star\in \text{SO}(3)$ for a rigid body. The kinematics of the rigid body are given by

\vspace{-0.4cm}\noindent
\begin{align}\label{eq:rigid_body_system}
    \dot{R}&= \sum_{i=1}^3 R \,\widehat{e}_i \,\theta_i \, u_i, 
\end{align}

\vspace{-0.4cm}\noindent
where $R\in\text{SO}(3)\subset\mathbb{R}^{3\times 3}$ is the rotation matrix representing the attitude of the rigid body, $u=(u_1,u_2,u_3)\in\mathbb{R}^3$ are the control inputs, and $\theta_i$ are the unknown control gains. Since $\text{SO}(3)$ is a Lie Group, the problem of stabilizing any specific attitude $R^\star\in \text{SO}(3)$ is equivalent to stabilizing the identity element, i.e. $I\in\text{SO}(3)$. Therefore, without loss of generality, we only consider the case when $R^\star=I$. We remark that, herein, we consider $\text{SO}(3)$ as an embedded submanifold of $\mathbb{R}^{3\times 3}$ equipped with its Riemannian metric, i.e. the Frobenius norm. 

Following \cite{Mayhew10Thesis}, we introduce the family of functions
\begin{align}
    \tilde{V}_q(R) &:= W\circ \Phi_q(R), & W(R)&:= \text{tr}(A (I-R)),
\end{align}
where $q\in\mathcal{Q}=\{1,2\}$, $A$ is the matrix given by
\begin{align}
    A = \frac{3}{\sum_{i=1}^3 \langle \tilde{\omega},e_i\rangle} \sum_{i=1}^3 \langle \tilde{\omega},e_i\rangle e_i,
\end{align}
$\tilde{\omega} = (11,12,13)$, the maps $\Phi_q:\text{SO}(3)\rightarrow\text{SO}(3)$ are
\begin{align}
    \Phi_q(R):=\exp\left(\frac{(3-2q)}{4} W(R)[\omega]_{\times}\right)R,
\end{align} 
and $\omega = \tilde{\omega}/\lVert \tilde{\omega}\lVert \in\mathbb{S}^2$. As shown in \cite{Mayhew10Thesis}, the family of functions $\{\tilde{V}_1,\tilde{V}_2\}$ is a $\delta$-gap synergistic family of potential functions with respect to $I$ for any $\delta\in(0,1/2)$. 

Next, let $\mathcal{M}$ be the Euclidean submanifold
\begin{align}
    \mathcal{M}=\{{vec}(R)~|~ R\in\text{SO}(3)\} \subset\mathbb{R}^9.
\end{align}
By defining $p={vec}(R)\in\mathcal{M}$, it follows that $p$ evolves according to the driftless control-affine system
\begin{align}
    \dot{p}&= \sum_{i=1}^3 b_i(p)\,\theta_i\,u_i, & b_i(p) &= -(\widehat{e}_i\otimes I)p
\end{align}
Since the map ${vec}:\text{SO}(3)\rightarrow\mathcal{M}$ is a diffeomorphism, the family of functions $\{\tilde{V}_1,\tilde{V}_2\}$ can be pulled back (via the inverse of $vec$) to a $\delta$-gap synergistic family of potential functions with respect to $p^\star={vec}(I)$. More explicitly, the family of functions $\{V_1,V_2\}$ defined by $V_q(p) = \tilde{V}_q\circ vec^{-1}(p)$ is a $\delta$-gap synergistic family of potential functions with respect to $p^\star={vec}(I)$. Therefore, it follows from Proposition \ref{prop:SCLF_synergistic} that the function $V$ defined by $V(p,q)=V_q(p)$ is an SCLF with respect to $\mathcal{A}=\{p^\star\}\times\mathcal{Q}$ for the HDS $\mathcal{H}$, and therefore that $\mathcal{A}\times\Theta\times\mathbb{T}^3$ is SGpAS for $\mathcal{H}_{cl}$.

We now provide numerical simulations. The results of the simulations are shown in Figure \ref{fig:rigid_body}. To generate the results, we used $\gamma = \sqrt{1}$, $\kappa = 4$, $T_1=1$, $T_2=2$, $T_3=3$, $\delta = 1/5$, and $\varepsilon=1/\sqrt{12\pi}$. The target point is $p^\star=vec(I)=(1,0,0,0,1,0,0,0,1)$. To emphasize the robustness of the proposed feedback law compared to the non-hybrid model-free controllers \cite{DurrManifold}, we added a small adversarial perturbation that locally stabilizes the (bad) critical point $p^\sharp=(-1,0,0,0,1,0,0,0,-1)$ in the absence of switching. As shown in Figure \ref{fig:rigid_body}, the proposed control law is not affected by the perturbation whereas the vanilla model-free feedback law is trapped by the adversarial perturbation in the vicinity of the critical point $p^\sharp$.

\subsection{Robust Target Seeking with Obstacle Avoidance for a Nonholonomic Vehicle}\label{subsec:obstacle_avoidance_nonholonomic}
Finally, we consider again the problem of robust global stabilization of a target position for a mobile vehicle moving in an obstructed planar domain. However, in contrast to the single integrator dynamics considered in Section 4.2, we now consider a \emph{nonholonomic} vehicle model. Due to the nonholonomic nature of the vehicle, the control vector fields do not span the entire tangent space of the state space manifold, and therefore Proposition \ref{prop:SCLF_synergistic} is not applicable. Although we are still able to apply the model-free feedback law \eqref{eq:ES_Law}, the stability analysis is different from the examples considered hitherto. Thus, we treat this special case separately.

Let $\mathcal{O}\subset\mathbb{R}^2$ denote the obstacle and $z=(z_1,z_2)\in\mathbb{R}^2\backslash\mathcal{O}$ denote the position of the vehicle. Suppose the vehicle is governed by the \emph{nonholonomic} kinematic equations
\begin{align}\label{nonholonomicmodel}
    \dot{z}&= \theta_1 u_1\psi, & \dot{\psi}&= u_2 S \psi
\end{align}
where $u=(u_1,u_2)\in\mathbb{R}^2$ are the control inputs, $\theta_1\in \mathcal{E}\subset\mathbb{R}$ is the unknown control gain, $S$ is the matrix defined in \eqref{eq:single_oscillator_dynamics}, and $\psi\in\mathbb{S}^1$ is a unit vector that indicates the forward direction of motion, i.e. the orientation, of the vehicle. This kinematic model has also been considered in \cite{CochranSourceSeekingTAC09,Poveda:20TAC,zhang2007source,yang1999sliding,durrEbenbauer}. The goal of the vehicle is to stabilize a known target position $z^\star=(z_1^\star,z_2^\star)\in\mathbb{R}^2\backslash\mathcal{O}$, irrespective of the orientation $\psi\in\mathbb{S}^1$, while avoiding the obstacle $\mathcal{O}$. We impose Assumption \ref{asmp:obstacle_contained} on the position of the target point $z^\star$ relative to the obstacle.

\begin{figure}[t]
    \centering
    \includegraphics[width=\linewidth]{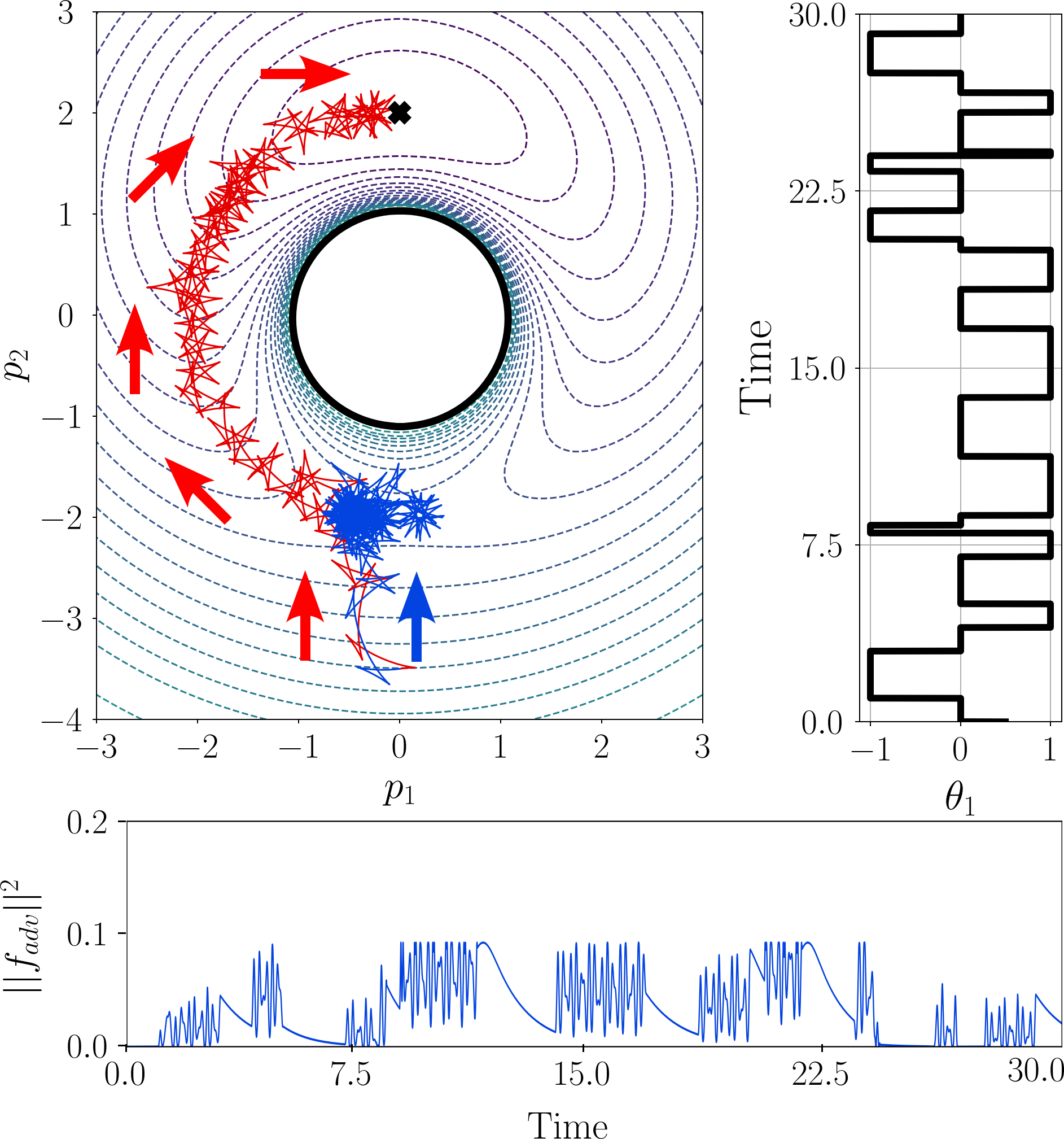}
    \caption{Simulation results for the problem of target-stabilization with obstacle avoidance for the \emph{nonholonomic kinematics} \eqref{eq:obstacle_avoidance_kinematics_nonholonomic}. }
    \label{fig:obstacle_avoidance}
\end{figure}

Under the diffeomorphism defined by \eqref{eq:obstacle_avoidance_diffeomorphism}, the pushforward of the nonholonomic kinematics of the vehicle is given by
\begin{align}\label{eq:obstacle_avoidance_kinematics_nonholonomic}
    \dot{p}&= \theta_1 u_1\text{D}\varphi\circ\varphi^{-1}(p)\psi, & \dot{\psi}&= u_2 S \psi,
\end{align}
where $p=(\rho,\vartheta)\in\mathbb{R}\times\mathbb{S}^1$. Therefore, the goal of globally stabilizing the subset $\{z^\star\}\times\mathbb{S}^1$ in $\big(\mathbb{R}^2\backslash(z_{\mathcal{O}}+d_0\,\mathbb{B}))\big)\times\mathbb{S}^1$ is equivalent to globally stabilizing the compact subset $\{p^\star\}\times\mathbb{S}^1$ on the smooth manifold $\varphi\big(\mathbb{R}^2\backslash(z_{\mathcal{O}}+d_0\,\mathbb{B}))\big)\times\mathbb{S}^1=\mathcal{M}\times\mathbb{S}^1$, where $p^\star=(\rho^\star,\vartheta^\star)=\varphi(z^\star)$. In contrast to the single integrator kinematics considered in subsection \ref{subsec:obstacle_avoidance_integrator}, Proposition \ref{prop:SCLF_synergistic} is not applicable for a vehicle with the nonholonomic kinematics \eqref{eq:obstacle_avoidance_kinematics_nonholonomic} since the control vector fields in \eqref{eq:obstacle_avoidance_kinematics_nonholonomic} clearly do not span the entire tangent space of the manifold $\mathcal{M}$ and so Assumption \ref{asmp:synergistic_manifold_regularity} is violated. Nevertheless, we will show that the structure of the control-affine system \eqref{eq:obstacle_avoidance_kinematics_nonholonomic} permits the use of a slightly different model-free feedback law that stabilizes $\mathcal{H}$. To that end, we introduce the family of functions
\begin{align*}
    V_q(\mu)&:= \frac{1}{2}(\rho-\rho^\star)^2 + \sqrt{(\text{e}^{\rho}-\text{e}^{\rho^\star})^2+1}-1 + W_q(\vartheta),
\end{align*}
where $q\in\mathcal{Q}=\{1,2\}$, and $W_q$ are the functions defined in \eqref{eq:syn_family_circle}. As in Section \ref{subsec:obstacle_avoidance_integrator}, the family of functions $\{V_1,V_2\}$ is a $\delta$-gap synergistic family of potential functions with respect to $\mu^\star=(\rho^\star,\vartheta^\star)$ on $\mathbb{R}\times\mathbb{S}^1$. Let 
$$\mu:=(\rho,\vartheta)\in\mathbb{R}\times\mathbb{S}^1,$$ 
$p=(\mu,\psi)\in\mathbb{R}\times\mathbb{S}^1\times\mathbb{S}^1 = \mathcal{M}$, $ x=(p,q)\in\mathcal{M}\times\mathcal{Q}$, and define the function $V:\mathcal{M}\times\mathcal{Q}\rightarrow\mathbb{R}_{\geq 0}$ as follows:
\begin{equation}\label{secondefinitionVcontrollaw}
V(x)=V_q(\mu).
\end{equation}
Also, let $r=1$, $T_1=1$, $\eta = \eta_1\in\mathbb{S}^1$, and let the first control input $u_1$ be given by the feedback law \eqref{eq:ES_Law}, i.e. we take $u_1$ to be
\begin{subequations}\label{eq:obstacle_avoidance_kinematics_nonholonomic_ES_law}
    \begin{align}
    u^\varepsilon_1(V(x),\eta) = \varepsilon^{-1}\sqrt{\frac{4\pi\gamma}{\kappa}} \langle\exp(\kappa V(x) S)e_1,\eta_1 \rangle.
    \end{align}
    Next, let the second control input $u_2$ to be defined as:
    \begin{align}
        u^\varepsilon_2 = 2\pi\varepsilon^{-1}.
    \end{align}
\end{subequations}
Finally, define the compact sets $\mathcal{A}=\{\mu^\star\}\times\mathbb{S}^1\times\mathcal{Q}$, and $\bar{\mathcal{A}}=\mathcal{A}\times
\Theta$. Then, we have the following proposition proved in Section \ref{sec:Proofs}.
\begin{prop}\label{prop:obstacle_avoidance_nonholonomic_SCLF}\normalfont
    $\bar{\mathcal{A}}\times\mathbb{S}^1$ is SGpAS as $\varepsilon\rightarrow 0^+$ for $\mathcal{H}_{cl}$ with $V$ given by \eqref{secondefinitionVcontrollaw}. 
\end{prop}
%


We conclude this section with a numerical simulation result illustrating the control of the non-holonomic vehicle \eqref{nonholonomicmodel}. The resulting trajectories of the vehicle are shown in Figure \ref{fig:obstacle_avoidance}. To generate the results, we used $\gamma = 2$, $\delta = 1/4$, $\kappa = 4$, and $\varepsilon=1/\sqrt{6\pi}\approx 0.165$. The target position is $z^\star=(0,2)$, and we considered an obstacle with radius $d=1$, centered at the origin, i.e. $z_{\mathcal{O}}=(0,0)$. The system is simulated under an additive adversarial perturbation designed to trap the trajectories of the vehicle behind the obstacle. The red trajectory corresponds to the proposed hybrid model-free feedback law, whereas the blue trajectory corresponds to a non-hybrid, model-free feedback law \cite{scheinker2017model,durrEbenbauer}. The arrows indicate the direction of motion for the corresponding trajectory. The figure on the bottom depicts the magnitude of the adversarial perturbation which effectively traps the vanilla model-free feedback behind the obstacle. The figure on the top right depicts the control gain $\theta_1$ as a function of time. The figure on the top left depicts the obstructed planar domain wherein the vehicle operates. 

\section{Proofs}
\label{sec:Proofs}
In this section, we present the proofs of our main results.
\subsection{Proof of Theorem \ref{thm:main_theorem}}
Let $\tau\in\mathbb{R}_{\geq 0}$ and consider the HDS $\widetilde{\mathcal{H}}_{cl}$ with state $(\xi,\eta,\tau)\in\mathbb{R}^n\times\mathbb{R}^r\times\mathbb{R}_{\geq0}$ and dynamics:
\begin{align*}
    \begin{cases}
                C\times\mathbb{T}^r\times\mathbb{R}_{\geq 0},&
                \begin{pmatrix}
                    \dot{\xi}\\\dot{\eta}\\\dot{\tau}
                \end{pmatrix}\in F_\varepsilon(\xi,\eta)\times \{\Lambda_\varepsilon(\eta)\}\times\{\varepsilon^{-2}\}
            \\
            D\times\mathbb{T}^r\times\mathbb{R}_{\geq 0}, & \begin{pmatrix}
                \xi^+\\
                \eta^+\\
                \tau^+
            \end{pmatrix}\in G(\xi)\times\{\eta\}\times\{\tau\}
    \end{cases}
\end{align*}
where $F_{\varepsilon}$ is given by \eqref{fsubespilon}, $\Lambda_\varepsilon(\eta)$ is given by \eqref{eq:lambda}, $G$ is given by \eqref{jumpmapGd}, and $C,D$ are given by \eqref{constructionflowjumpsets}. System $\widetilde{\mathcal{H}}_{cl}$ is a trivial dynamic extension of the closed-loop HDS $\mathcal{H}_{cl}$. Therefore, given an initial condition, any solution of the HDS $\mathcal{H}_{cl}$ corresponds to some solution of the HDS $\widetilde{\mathcal{H}}_{cl}$. Let $\zeta:=(\zeta_1,\zeta_2,\dots,\zeta_r)\in\mathbb{S}^1\times\mathbb{S}^1\times\cdots\mathbb{S}^1=\mathbb{T}^r$ and consider the HDS $\widehat{\mathcal{H}}_{cl}$ defined by
\begin{align}\label{auxiliarysystem001}
        \begin{cases}
                C\times\mathbb{T}^r\times\mathbb{R}_{\geq 0},&
                \begin{pmatrix}
                    \dot{\xi}\\\dot{\zeta}\\\dot{\tau}
                \end{pmatrix}\in \tilde{F}_\varepsilon(\xi,\zeta,\tau)\times \{0\}\times\{\varepsilon^{-2}\}
            \\
            D\times\mathbb{T}^r\times\mathbb{R}_{\geq 0}, & \begin{pmatrix}
                \xi^+\\
                \zeta^+\\
                \tau^+
            \end{pmatrix}\in G(\xi)\times\{\zeta\}\times\{\tau\}
            \end{cases}
\end{align}

\vspace{-0.5cm}
where $\tilde{F}_\varepsilon$ is defined by
\begin{align*}
        \tilde{F}_\varepsilon(\xi,\zeta,\tau):= F_{\varepsilon}(\xi,\exp( \Omega\tau)\zeta),
\end{align*}
and the matrix $\Omega$ is a block diagonal matrix of size $2r\times 2r$ with $r$ diagonal blocks of size $2\times 2$ such that the $i$-th block is the matrix $2\pi T_i^{-1} S$, where $S$ is the matrix defined in \eqref{eq:ES_Law}. Let $(\xi,\eta,\tau):\text{dom}(\xi,\eta,\tau)\rightarrow \mathbb{R}^n\times\mathbb{T}^r\times\mathbb{R}_{\geq 0}$ be any solution of $\widetilde{\mathcal{H}}_{cl}$ and define the hybrid arc $(\xi,\zeta,\tau):\text{dom}(\xi,\eta,\tau)\rightarrow \mathbb{R}^n\times\mathbb{T}^r\times\mathbb{R}_{\geq 0}$ by
    \begin{align}\label{equivalenceofsolutions}
        (\xi,\zeta,\tau):= (\xi,\exp(-\Omega\tau)\eta,\tau).
    \end{align}
    Direct computation shows that
    \begin{align}\label{equivalenceofsolutions}
        \dot{\zeta}&=\exp(-\Omega\tau)\dot{\eta} - \dot{\tau}\exp(-\Omega\tau) \Omega\eta = 0.
    \end{align}
    Since, for any $\tau\in\mathbb{R}_{\geq 0}$ and any $\eta$ in $\mathbb{T}^r$, $\exp(-\Omega\tau)\eta \in \mathbb{T}^r$, it follows that the hybrid arc $(\xi,\zeta,\tau)$ defined by \eqref{equivalenceofsolutions} is a solution of $\widehat{\mathcal{H}}_{cl}$ with the same initial condition for $\xi$ and $\tau$, and the following initial condition for $\zeta$:
    \begin{align}
        \zeta(0,0) = \exp(-\Omega\tau(0,0))\eta(0,0).
    \end{align}
    Hence, for any solution of the HDS $\widetilde{\mathcal{H}}_{cl}$ there exists some solution of the HDS $\widehat{\mathcal{H}}_{cl}$ such that \eqref{equivalenceofsolutions} holds. In particular, the $\xi$-component of both solutions coincide. We will use this equivalence to establish that every solution of $\widetilde{\mathcal{H}}_{cl}$ satisfies suitable stability bounds by first establishing such bounds for \emph{every} solution of the HDS $\widehat{\mathcal{H}}_{cl}$.
    
    An explicit computation of $\tilde{F}_\varepsilon$ shows that
    \begin{align}
        \tilde{F}_\varepsilon(\xi,\zeta,\tau) = \{\tilde{f}_\varepsilon(\xi,\zeta,\tau)\}\times F_e(\theta),
    \end{align}
    where $F_e$ is given by \eqref{eq:exogenous_signal_dynamics}, and
    \begin{align}
        \tilde{f}_\varepsilon(\xi,\zeta,\tau) = f_0(x,\theta) + \sum_{i=1}^r f_i(x,\theta) \hat{u}_i(V(x),\zeta_i,\tau),
    \end{align}
    and the functions $\hat{u}_i$ are given by
    \begin{align*}
        \hat{u}_i(V,\zeta_i,\tau) &= \varepsilon^{-1}\sqrt{\frac{4\pi\gamma}{T_i\kappa}} \langle \exp(\kappa V S)e_1,\exp(2\pi S \tau/T_i)\zeta_i\rangle. \notag\\
    \end{align*}
    Now, since the constants $\{T_1,T_2,\dots,T_r\}$ are rational, it follows that there exists $T\in\mathbb{R}_{>0}$ such that the map $\tilde{f}_\varepsilon$ is $T$-periodic in $\tau$. Moreover, by definition of the map $\tilde{f}_{\varepsilon}$, the HDS $\widehat{\mathcal{H}}_{cl}$ belongs to the class of well-posed HDS with highly oscillatory flow maps which was analyzed in \cite{abdelgalil2023lie}. By applying the Lie-bracket averaging results in \cite{abdelgalil2023lie}, we obtain that the Hybrid Lie-bracket averaged system $\widehat{\mathcal{H}}_{cl}^{\text{ave}}$ corresponding to the HDS $\widehat{\mathcal{H}}_{cl}$ is given 
    by
    %
    \begin{align}\label{averagehybridsystem0}
        \widehat{\mathcal{H}}_{cl}^\text{ave}:\begin{cases}
            ~C\times \mathbb{T}^r,&
            (\dot{\xi},\dot{\zeta})\hphantom{^+}\in \bar{F}(\xi)\times\{0\}
        \\
        ~D \times \mathbb{T}^r, & (\xi,\zeta)^+\in G(\xi)\times\{\zeta\}
        \end{cases},
    \end{align}
    wherein the flow map $\bar{F}$ is defined by
    \begin{align*}
        \bar{F}(\xi)&:= \{\bar{f}(x,\theta)\}\times F_e(\theta), \\
        \bar{f}(\xi)&:= f_0(x,\theta) - \gamma {\sum_{i=1}^r}\langle\nabla V(x),f_i(x,\theta)\rangle f_i(x,\theta),
    \end{align*}
    %
    which is independent of $\zeta$. Thus, the HDS $\widehat{\mathcal{H}}_{cl}^{\text{ave}}$ is nothing but a trivial dynamic extension of the HDS $\mathcal{H}_V$ obtained by adding $\zeta$ as a state with trivial flow and jump dynamics, i.e., $\dot{\zeta}=0$ and $\zeta^+=\zeta$. From the assumptions of the theorem, namely that $\mathcal{H}$ is strongly $\nabla V$-stabilizable, it follows that the subset $\mathcal{A}\times(\Theta_C\cup\Theta_D)$ is UGAS for the HDS $\mathcal{H}_V$ defined in \eqref{eq:org_HDS_V}. It follows that the compact subset $\bar{\mathcal{A}}$ is UGAS for the HDS $\widehat{\mathcal{H}}_{cl}^{\text{ave}}$ given by \eqref{averagehybridsystem0}. From \cite[Theorem 2]{abdelgalil2023lie}, we obtain that there exists a class $\mathcal{KL}$ function $\beta$ such that for each compact subset $K\subset (C\cup D)\times\mathbb{T}^{r}$ and for each $\nu > 0$, there exists $\varepsilon^*>0$ such that for all $\varepsilon\in(0,\varepsilon^*]$ and for all solutions of $\widehat{\mathcal{H}}_{cl}$ with $(\xi(0,0),\zeta(0,0))\in K$, the following inequality holds for all $(t,j)\in\text{dom}(x,\zeta,\tau)$:
    \begin{align*}
        |(\xi(t,j),\zeta(t,j))|_{\bar{\mathcal{A}}}\leq \beta(|(\xi(0,0),\zeta(0,0))|_{\bar{\mathcal{A}}},t+j)+ \nu.
    \end{align*}    
    With $\varepsilon\in (0,\varepsilon^*]$, let $(\xi,\eta,\tau):\text{dom}(\xi,\eta,\tau)\rightarrow \mathbb{R}^n\times\mathbb{T}^r\times\mathbb{R}_{\geq 0}$ be any solution of the HDS $\widetilde{\mathcal{H}}_{cl}$ such that $(\xi(0,0),\eta(0,0))\in K$. By construction, it follows that $(\xi,\eta,\tau)$ satisfies 
    \begin{align}
        (\xi,\eta,\tau)= (\xi,\exp(\Omega\tau)\zeta,\tau),
    \end{align}
    for some solution $(\xi,\zeta,\tau)$ of the HDS $\widehat{\mathcal{H}}_{cl}$ with $(\xi,\zeta)(0,0)\in K$. Therefore, we obtain that there exists a class $\mathcal{KL}$ function $\beta$ such that for each compact subset $K\subset (C\cup D)\times\mathbb{T}^{r}$ and for each $\nu > 0$, there exists $\varepsilon^*>0$ such that for all $\varepsilon\in(0,\varepsilon^*]$ and for all solutions of $\widetilde{\mathcal{H}}_{cl}$ with $(\xi(0,0),\eta(0,0))\in K$, the following inequality holds for all $(t,j)\in\text{dom}(x,\eta,\tau)$:
    \begin{align*}
        |(\xi(t,j),\eta(t,j))|_{\bar{\mathcal{A}}}\leq \beta(|(\xi(0,0),\eta(0,0))|_{\bar{\mathcal{A}}},t+j)+ \nu.
    \end{align*} 
    However, since $\widetilde{\mathcal{H}}_{cl}$ is nothing but a trivial dynamic extension of the closed loop system $\mathcal{H}_{cl}$, it follows that the $\mathcal{KL}$ bound above is also true for every solution of $\mathcal{H}_{cl}$ with $(\xi(0,0),\eta(0,0))\in K$. In particular, we may always take $K= K_1\times\mathbb{T}^r$ for some compact $K_1\subset (C\cup D)$ so that any initial condition $\eta(0,0)$ for the state of the oscillators is admissible. This concludes the proof.
    \subsection{Proof of Proposition \ref{prop:SCLF_synergistic}}
    When $f_0$ and $f_i$ are as given in \eqref{eq:synergistic_mnfd_data}, direct computations show that $\bar{f}$, defined in \eqref{eq:avg_control_affine_sys}, takes the form
    \begin{align*}
        \bar{f}(\xi) = -\gamma\sum_{i=1}^r\theta_i^2\langle\nabla V_q(p),b_i(p)\rangle\begin{pmatrix} b_i(p)\\0 \end{pmatrix}.
    \end{align*}
     As a result, if $\xi=(x,\theta) = ((p,q),\theta)$, we have that
    \begin{align*}
        \dot{V}(\xi) \leq -\gamma \sum_{i=1}^r \theta_i^2 \langle \nabla V_q(p),b_i(p)\rangle ^2 \leq 0,
    \end{align*}
    for all $\xi\in C$. On the other hand, by construction we have that for all $\xi\in D$:
    \begin{align*}
        \Delta V(\xi)\leq 0,
    \end{align*}
    Thus, $V$ is an SCLF candidate with respect to $\mathcal{A}$ for the HDS $\mathcal{H}$ defined by \eqref{eq:org_HDS}-\eqref{eq:combined_org_HDS}, \eqref{hybridautomaton}, and \eqref{eq:synergistic_mnfd_data}. 
    
    Since $\mathcal{M}$ is equipped with the Riemannian metric of the ambient space and $b_i(p)\in T_p\mathcal{M}$, we have that
    \begin{align*}
        \langle \nabla V_q(p),b_i(p)\rangle = \langle \nabla_{\mathcal{M}} V_q(p),b_i(p)\rangle,
    \end{align*}
    where $\nabla_{\mathcal{M}} V_q(p)$ is the unique orthogonal projection of $\nabla V_q(p)$ onto the tangent space $T_p\mathcal{M}$ with respect to the Riemannian metric of the ambient Euclidean space. Therefore, we obtain that
    \begin{align*}
        \dot{V}(\xi) = -\gamma \sum_{i=1}^r \theta_i^2\langle \nabla_{\mathcal{M}} V_q(p),b_i(p)\rangle^2.
    \end{align*}
    On the other hand, by definition, if $\theta\in(\Theta_C\cup\Theta_D)\backslash \mathcal{E}_b$, then $\theta_i^2> 0$, $\forall i\in\{1,2,\dots,r\}$. Since $\mathcal{E}$ is a discrete finite set, there exists a constant $\lambda_1\in\mathbb{R}_{>0}$ such that, for all $\theta\in(\Theta_C\cup\Theta_D)\backslash \mathcal{E}_b$, we have $\theta_i^2> \lambda_1$  $\forall i\in\{1,2,\dots,r\}$.
 
    Hence, for all $\xi\in C$ such that $\theta\in(\Theta_C\cup\Theta_D)\backslash \mathcal{E}_b$, we have $\dot{V}(\xi) \leq -\gamma\lambda_1 \sum_{i=1}^r \langle \nabla_{\mathcal{M}} V_q(p),b_i(p)\rangle^2$. From Assumption \ref{asmp:synergistic_manifold_regularity}, there exists a constant $\lambda\in\mathbb{R}_{>0}$ such that, for all $p\in\mathcal{M}$ and all $v\in T_p\mathcal{M}$, we have that $\sum_{i=1}^r \langle v,b_i(p)\rangle^2 \geq \lambda\langle v,v\rangle$. Hence, for all $\xi\in C$ such that $\theta\in(\Theta_C\cup\Theta_D)\backslash \mathcal{E}_b$, we have
    \begin{align*}
        \dot{V}(\xi) \leq -\gamma\lambda_1\lambda \langle \nabla_{\mathcal{M}} V_q(p),\nabla_{\mathcal{M}} V_q(p)\rangle .
    \end{align*}
    By construction of the family $\{V_q\}_{q\in\mathcal{Q}}$, we have that
    \begin{align*}
        \nabla_{\mathcal{M}} V_q(p)=0, \iff p = p^\star,
    \end{align*}
    for all $(p,q)\in \mathcal{X}_C\times\mathcal{Q}$. It follows that there exists a positive definite function $\rho$ such that
    \begin{align*}
        \langle \nabla_{\mathcal{M}} V_q(p),\nabla_{\mathcal{M}} V_q(p)\rangle \geq \rho(|(p,q)|_{\mathcal{A}}),
    \end{align*}
    for all $(p,q)\in \mathcal{X}_C$, which implies that for all $\xi\in C$ such that $\theta\in(\Theta_C\cup\Theta_D)\backslash \mathcal{E}_b$, we have that
    \begin{align*}
        \dot{V}(\xi) \leq -\gamma\lambda_1\lambda \rho(|x|_{\mathcal{A}}) .
    \end{align*}
    Following similar steps as in \cite{kellett2014compendium}, it can be shown that there exists a continuously differentiable $\mathcal{K}_\infty$ function $\bar{\alpha}$ such that the function $\hat{V} = \bar{\alpha}\circ V$ satisfies
    \begin{align*}
        \bar{\alpha}_1(|x|_{\mathcal{A}})\leq &\hat{V}(x) \leq \bar{\alpha}_2(|x|_{\mathcal{A}}), & \forall x&\in \mathcal{X}, \\
        &\dot{\hat{V}}(\xi)\leq 0, & \forall \xi &\in C,\\
        \Delta &\hat{V}(\xi)\leq 0, & \forall \xi&\in D.
    \end{align*}
    and, for all $\xi\in C$ such that $\theta\in (\Theta_C\cup\Theta_D)\backslash \mathcal{E}_b$, we have
    \begin{align*}
        \dot{\hat{V}}(\xi) \leq - \hat{V}(\xi).
    \end{align*}
    We now introduce the function $\tilde{V}(\xi)= \hat{V}(\xi) \text{e}^{\theta_{r+2}}$. Direct computations shows that, during flows, we have
    \begin{align*}
        \dot{\tilde{V}}(\xi) &= (\dot{\hat{V}}(\xi) + \dot{\theta}_{r+2} \hat{V}(\xi))\text{e}^{\theta_{r+2}}
    \end{align*}
    From \eqref{hybridautomaton}, we have that $\dot{\theta}_{r+2} \in [0,\chi_2]-\mathbb{I}_{\mathcal{E}_b}(\theta)$. Therefore, if $\theta\in\mathcal{E}_b$, we have that
    \begin{align*}
        \dot{\theta}_{r+2}\leq -(1-\chi_2) \implies \dot{\tilde{V}}(\xi) \leq -(1-\chi_2) \tilde{V}(\xi).
    \end{align*}
    On the other hand, if  $\theta\in (\Theta_C\cup\Theta_D)\backslash \mathcal{E}_b$, we also have that
    \begin{align*}
        \dot{\theta}_{r+2}\leq \chi_2 \implies \dot{\tilde{V}}(\xi) \leq -(1-\chi_2) \tilde{V}(\xi)
    \end{align*}
    From \eqref{hybridautomaton}, we have that $(\theta_{r+2})^+ = \theta_{r+2}$, for all $\xi\in D$. Thus, using the facts that $\theta_{r+2}\in[0,T_\circ]$ and $|x|_{\mathcal{A}} = |\xi|_{\tilde{\mathcal{A}}}$ for all $\xi\in C\cup D$, we obtain 
    \begin{subequations}\label{eq:modified_SCLF_proof}
        \begin{align}
            \bar{\alpha}_1(|\xi|_{\tilde{\mathcal{A}}})\leq &\tilde{V}(\xi) \leq \bar{\alpha}_2(|\xi|_{\tilde{\mathcal{A}}})\text{e}^{T_\circ}, & \forall \xi&\in C\cup D, \\
            \dot{\tilde{V}}(\xi)\leq &-(1-\chi_2)\tilde{V}(\xi), & \forall \xi &\in C,\\
            \Delta &\tilde{V}(\xi)\leq 0, & \forall \xi&\in D,
        \end{align}
    \end{subequations}

    \vspace{-0.4cm}
    \noindent 
    where $\tilde{\mathcal{A}} = \mathcal{A}\times(\Theta_C\cup\Theta_D)$. 
    Because the set $\tilde{\mathcal{A}}$ is Lyapunov stable for $\mathcal{H}_V$, which is entailed by the conditions \eqref{eq:modified_SCLF_proof} \cite[Theorem 3.18]{RSanfeliceBook}, it follows that any maximal solution to $\mathcal{H}_V$ starting in $C\cup D$ is bounded. Thanks to Lemma \ref{lem:solution_properties_lemma_aux}, we obtain that, by invoking \cite[Proposition 2.34]{RSanfeliceBook}, no nontrivial maximal solution to $\mathcal{H}_V$ starting in $C\cup D$ is Zeno and that, for every maximal solution $\xi$ to $\mathcal{H}_V$ starting in $C\cup D$, there exists $t_\circ>0$ such that $\overline{t}_{j}-\underline{t}_{j}\geq t_\circ > 0$, for all $j\geq 1$, where $\overline{t}_{j}:=\sup\{t\in\mathbb{R}_{\geq 0}~|~(t,j)\in\text{dom}(\xi)\}$, and
    $\underline{t}_{j}:=\,\inf\,\{t\in\mathbb{R}_{\geq 0}~|~(t,j)\in\text{dom}(\xi)\}$. In other words, no solution of $\mathcal{H}_V$ is discrete. By invoking \cite[Theorem 3.19]{RSanfeliceBook}, we obtain that the compact subset $\tilde{\mathcal{A}}$ is globally asymptotically stable in the sense of \cite[Definition 3.1]{RSanfeliceBook}. However, since $\tilde{\mathcal{A}}$ is compact and $\mathcal{H}_V$ is a well-posed HDS, it follows from \cite[Theorem 3.22]{RSanfeliceBook} that the global asymptotic stability of $\tilde{\mathcal{A}}$ is \emph{uniform}, i.e. that $\tilde{\mathcal{A}}$ is UGAS for $\mathcal{H}_V$ in the sense of Definition \ref{definitionstablity1} which concludes the proof.

    \vspace{-0.2cm}
    \subsection{Proof of Proposition \ref{prop:obstacle_avoidance_nonholonomic_SCLF}}
    
    \vspace{-0.2cm}
    The closed-loop HDS $\mathcal{H}_{cl}$ is defined by \eqref{eq:combined_org_HDS} in conjunction with \eqref{eq:exogenous_signal_dynamics}, \eqref{eq:ES_Law}, \eqref{eq:closed_loop_org_HDS}, \eqref{hybridautomaton}, \eqref{eq:synergistic_mnfd_data}, \eqref{eq:obstacle_avoidance_kinematics_nonholonomic}, and \eqref{eq:obstacle_avoidance_kinematics_nonholonomic_ES_law}.
    Let $\tau_1,\tau_2\in\mathbb{R}_{\geq 0}$, $\tau=(\tau_1,\tau_2)\in\mathbb{R}_{\geq 0}^2 = \mathbb{R}_{\geq 0}\times\mathbb{R}_{\geq 0}$, and define the HDS $\widetilde{\mathcal{H}}_{cl}$ via the following data:
    \begin{align*}
        \begin{cases}
                ~\tilde{C}\times\mathbb{R}^2_{\geq 0},&
                \begin{pmatrix}
                    \dot{\xi}\\\dot{\eta}\\\dot{\tau}
                \end{pmatrix}\in F_\varepsilon(\xi,\eta)\times \{\Lambda_\varepsilon(\eta)\}\times\left\{(\varepsilon^{-1},\varepsilon^{-2})\right\}
            \\
            ~\tilde{D}\times\mathbb{R}_{\geq 0}^2, & \begin{pmatrix}
                \xi^+\\
                \eta^+\\
                \tau^+
            \end{pmatrix}\in G(\xi)\times\{\eta\}\times\{\tau\}
            \end{cases}
    \end{align*}
    where $\tilde{C} = C\times\mathbb{T}^1$, $\tilde{D}=D\times\mathbb{T}^1$. It is clear that $\widetilde{\mathcal{H}}_{cl}$ is a trivial dynamic extension of $\mathcal{H}_{cl}$. Next, similar to the proof of Theorem \ref{thm:main_theorem}, we introduce the HDS $\widehat{\mathcal{H}}_{cl}$ defined by
    \begin{align*}
        \begin{cases}
                \tilde{C}\times\mathbb{R}^2_{\geq 0}&
                \begin{pmatrix}
                    \dot{\mu}\\\dot{\zeta}\\\dot{\tau}
                \end{pmatrix}\in \widehat{F}_\varepsilon(\mu,\zeta,\tau)\times \{0\}\times\left\{(\varepsilon^{-1},\varepsilon^{-2})\right\}
            \\
            \tilde{D}\times\mathbb{R}_{\geq 0}^2 & \begin{pmatrix}
                \mu^+\\
                \zeta^+\\
                \tau^+
            \end{pmatrix}\in G(\mu)\times\{\zeta\}\times\{\tau\}
            \end{cases}
    \end{align*}
    where $\zeta\in\mathbb{S}^1$, $\mu=(p,\phi,q,\theta) \in \mathcal{M}\times\mathbb{S}^1\times\mathcal{Q}\times (\Theta_C\cup\Theta_D)$, and the map $\widehat{F}_\varepsilon$ is given by
    \begin{align*}
    \widehat{F}_\varepsilon(\mu,\zeta,\tau) = F_\varepsilon((p,\exp(2\pi\tau_1 S)\phi,q,\theta),\exp(2\pi \tau_2 S)\eta).
    \end{align*}
    An explicit computation shows that $\widehat{F}_\varepsilon(\mu,\zeta,\tau) = \{\widehat{f}_\varepsilon(\mu,\zeta,\tau)\}\times F_e(\theta)$, where the map $\widehat{f}_\varepsilon$ is defined by
    \begin{align}
        \widehat{f}_\varepsilon(\mu,\zeta,\tau) = \theta_1 b_1(p,\phi,\tau_1) \hat{u}_1^\varepsilon(V(x),\zeta,\tau_2),
    \end{align}
    where the vector field $b_1$ is given by
    \begin{align}
f_1(p,\phi,\tau_1):=\begin{pmatrix}\text{D}\varphi\circ\varphi^{-1}(p)\exp(2\pi\tau_1 S)\phi \\ 0 \\ 0\end{pmatrix},
    \end{align}
    and the function $\hat{u}_1^\varepsilon$ is given by
    \begin{align*}
        \hat{u}_1^\varepsilon(V,\zeta,\tau_2) &= \varepsilon^{-1}\sqrt{\frac{2\gamma}{\kappa}} \langle \exp(\kappa V S)e_1,\exp(2\pi S \tau_2)\zeta\rangle. \notag
    \end{align*}
    It is clear that $\widehat{f}_\varepsilon$ is $1$-periodic in $\tau_1$ and $\tau_2$. Moreover, the structure of the HDS $\widehat{\mathcal{H}}_{cl}$ is the same as the structure of the class of well-posed HDS with highly oscillatory flow maps considered in \cite{abdelgalil2023lie}. By applying the \emph{recursive} Lie-bracket averaging formulas in \cite{abdelgalil2023lie}, it can be shown that the Hybrid Lie-bracket averaged system $\widehat{\mathcal{H}}_{cl}^{\text{ave}}$ corresponding to $\widehat{\mathcal{H}}_{cl}$ is the system defined by
    \begin{align*}
        \begin{cases}
               (\mu,\zeta)\in \tilde{C},&
                 (\dot{\mu},\dot{\zeta})\hphantom{^+}\in \bar{F}(\mu)\times \{0\}
            \\
            (\mu,\zeta)\in\tilde{D}, & (\mu,\zeta)^+\in G(\mu)\times\{\zeta\},
            \end{cases}
    \end{align*}
    where the map $\bar{F}(\mu,\zeta)=\{\bar{f}(\mu)\}\times F_e(\theta)$,
    \begin{align}
        \bar{f}(\mu)&= -\frac{\theta_1^2\gamma}{2}\sum_{i=1}^2\langle b_i(p),\nabla V_q(p)\rangle \begin{pmatrix} b_i(p) \\ 0 \\ 0\end{pmatrix},
    \end{align}
    and the maps $b_i$ are given by $b_i(p)= \text{D}\varphi\circ\varphi^{-1}(p)e_i$. It is clear that $\widehat{\mathcal{H}}_{cl}^{\text{ave}}$ is a trivial dynamic extension of the HDS
    \begin{align*}
        \mathcal{H}_V: \begin{cases}
               ~\mu\in C,&
                 ~\dot{\mu}\hphantom{^+}\in \bar{F}(\mu)
            \\
           ~\mu\in  D, &~ \mu^+\in G(\mu)
            \end{cases}
    \end{align*}
    Now, consider the function $V(\mu) = V_q(p)$ and observe that
    \begin{align}
        \dot{V}(\mu) = -\frac{\theta_1^2\gamma}{2}\sum_{i=1}^2\langle b_i(p),\nabla V_q(p)\rangle^2 \leq 0,
    \end{align}
    for all $\mu\in D$, and   $\Delta V(\mu)\leq 0$, for all $\mu\in D$. In addition, $V$ is positive definite with respect to $\bar{\mathcal{A}} = \mathcal{A}\times (\Theta_C\cup\Theta_D)$ and, for every $c\in\mathbb{R}_{\geq 0}$, the set $\{\mu\in C\cup D~|~ V(\mu)\leq c\}$ is compact. Following similar steps to the proof of Proposition \ref{prop:SCLF_synergistic}, it can be shown that $\bar{\mathcal{A}}$ UGAS for the HDS $\mathcal{H}_V$, which implies that $\bar{\mathcal{A}}\times\mathbb{S}^1$ is UGAS for the HDS $\widehat{\mathcal{H}}_{cl}^{\text{ave}}$. Then, following similar steps as in the proof of Theorem \ref{thm:main_theorem}, we obtain that $\bar{\mathcal{A}}\times\mathbb{S}^1$ is SGpAS for $\mathcal{H}_{cl}$ as $\varepsilon\rightarrow 0^+$, which concludes the proof.
\section{Conclusions}
\label{sec:conclusions}
We studied the problem of robust global stabilization for a class of control-affine systems under dynamic uncertainty in the control directions and topological obstructions. By proposing a novel class of hybrid and highly-oscillatory feedback laws that seek the minimum of a family of synergistic Lyapunov functions, we provide a robust solution to this problem. Our method is particularly advantageous as it is model-free and can handle unknown and switching control gain signs, ensuring resilience in autonomous systems facing adversarial spoofing attacks. The practical relevance of our results is demonstrated through several concrete applications. Since our results allow to transfer stability properties from well-posed stable hybrid systems with \emph{known} control directions to systems with \emph{unknown} control directions, they open the door for the systematic development of new ``model-free'' controllers that can exploit the rich set of tools developed during the last two decades for hybrid systems. Numerical simulations validate the effectiveness and robustness of the proposed feedback law compared to existing techniques. Future research directions will focus on experimental validations of the proposed algorithms.

\vspace{-0.2cm}
\bibliographystyle{elsarticle-num}  
\bibliography{References}
\end{document}